\documentclass[twoside]{article}
\usepackage[accepted]{aistats2018}

\usepackage{amssymb}
\usepackage{mathrsfs}
\usepackage{amsfonts}
\usepackage{amsmath}
\usepackage{graphicx}
\usepackage{multirow}
\usepackage{caption2}
\usepackage{float}
\usepackage{booktabs}
\usepackage{algorithm}
\usepackage{algorithmic}
\usepackage{subfigure}
\usepackage{color}
\usepackage{tabu}

\newtheorem{Thm}{Theorem}
\newtheorem{Lemm}{Lemma}

\newtheorem{Coro}{Corollary}
\newtheorem{Assump}{Assumption}

\begin{document}



\twocolumn[

%\aistatstitle{Stochastic Nonconvex Semidefinite Optimization with Variance Reduction: Global Linear Convergence}
\aistatstitle{Finding Global Optima in Nonconvex Stochastic Semidefinite Optimization with Variance Reduction}

\aistatsauthor{ Jinshan Zeng$^{1,2}$ \And Ke Ma$^{3,4}$ \And  Yuan Yao$^{2,\dag}$ }

\aistatsaddress{
\textsuperscript{1} School of Computer Information Engineering, Jiangxi Normal University\\
\textsuperscript{2} Department of Mathematics, Hong Kong University of Science and Technology\\
\textsuperscript{3} State Key Laboratory of Information Security, Institute of Information Engineering, Chinese Academy of Sciences\\
\textsuperscript{4} School of Cyber Security, University of Chinese Academy of Sciences\\
\{jsh.zeng@gmail.com, make@iie.ac.cn, yuany@ust.hk\} ($\dag$ Corresponding author)
}
]

\begin{abstract}
There is a recent surge of interest in nonconvex reformulations via low-rank factorization for stochastic convex semidefinite optimization problem in the purpose of efficiency and scalability. Compared with the original convex formulations, the nonconvex ones typically involve much fewer variables, allowing them to scale to scenarios with millions of variables. However, it opens a new challenge that under what conditions the nonconvex stochastic algorithms may find the global optima effectively despite their empirical success in applications. In this paper, we provide an answer that a stochastic gradient descent method with variance reduction, can be adapted to solve the nonconvex reformulation of the original convex problem, with a \textit{global linear convergence}, i.e., converging to a global optimum exponentially fast, at a proper initial choice in the restricted strongly convex case. Experimental studies on both simulation and real-world applications on ordinal embedding are provided to show the effectiveness of the proposed algorithms.
\end{abstract}

\section{Introduction}

The stochastic convex semidefinite optimization problem, arising in many applications like non-metric multidimensional scaling \cite{Agarwal2007,Borg-NMDS2005}, matrix sensing \cite{Sanghavi2013,Recht-2016}, community detection \cite{Montanari-CommunityDetection2016}, synchronization \cite{Bandeira-Synchronization2016}, and phase retrieval \cite{Candes-ACHA2015}, is of the following form:
\begin{align}
\label{Eq:SD-stoch}
\mathop{\mathrm{min}}_{X \in \mathbb{R}^{p\times p}} \ f(X) = \frac{1}{n} \sum_{i=1}^n f_i(X)  \quad \mathrm{s.t.} \quad X \succeq 0,
\end{align}
where $f_i(X)$ is some convex, smooth cost function associated with the $i$-th sample, $X \succeq 0$ is the positive semidefinite (PSD) constraint.

There are many algorithms for solving problem \eqref{Eq:SD-stoch}, mainly including the first-order methods like the well-known projected gradient descent method \cite{Nestrov-1989}, interior point method \cite{Alizadeh-IPM1995}, and more specialized path-following interior point methods which use the (preconditioned) conjugate gradient or residual scheme to compute the Newton direction (for more detail, see the survey \cite{Monteiro-SDP2003} and references therein). However, most of these methods are not well-scalable due to the PSD constraint, i.e., $X \succeq 0$. To circumvent this difficult constraint, the idea of low-rank factorization was adopted in literature \cite{Monteiro2003,Monteiro2005} and became very popular in the past few years due to its empirical success \cite{Sanghavi-FGD2016}. Low-rank factorization recasts the original problem \eqref{Eq:SD-stoch} into an unconstrained problem by introducing another rectangular matrix $U\in \mathbb{R}^{p\times r}$ with $r < p$ such that $X = UU^T$. Let $g(U):= f(UU^T)$ and problem \eqref{Eq:SD-stoch} leads to,
\begin{align}
\label{Eq:SD-determ-ncvx}
\mathop{\mathrm{min}}_{U \in \mathbb{R}^{p\times r}} \ g(U)   \quad \mathrm{where} \quad r\leq p.
\end{align}

Problems \eqref{Eq:SD-stoch} and \eqref{Eq:SD-determ-ncvx} will be equivalent when $r\geq r^*$ in the sense that problem \eqref{Eq:SD-determ-ncvx} can find a global optimum $X^*$ of problem \eqref{Eq:SD-stoch} with $r^* = \mathrm{rank}(X^*)$.
Since the PSD constraint has been eliminated, the recast problem \eqref{Eq:SD-determ-ncvx} has a significant advantage over \eqref{Eq:SD-stoch}, but this benefit has a corresponding cost: the objective function is no longer convex but instead \textit{nonconvex} in general. Even for the simple first-order methods like the factored gradient descent (FGD), its \emph{global linear convergence}\footnote{By \textit{global linear convergence}, it means that the algorithm converges to a global optimum exponentially fast when the initial choice is in a prescribed ball.} remains unspecified until a recent work in \cite{Sanghavi-FGD2016}.
%Thus, in order to solve the original difficult problem \eqref{Eq:SD-stoch}, several first-order methods like  \cite{Sanghavi-FGD2016}, were proposed to solve the recast easier one \eqref{Eq:SD-determ-ncvx}.
%The linear \textit{global convergence}\footnote{The \textit{global convergence} used in this paper means the convergence to a global optimum.} of FGD was developed in \cite{Sanghavi-FGD2016}.
%Specifically,
%the Factored Gradient Descent (FGD) method can be described as follows: let $U^t$ be the iterate at the $t$-th iteration and $X^t:=U^t {U^t}^T$, then the next iterate $U^{t+1}$ is updated according to the following
%\begin{align}
%\label{alg:FGD-fixed}
%U^{t+1} = U^t - \eta \nabla f(X^t)U^t,
%\end{align}
%where $\eta >0$ is a step size.
%The \textit{global linear convergence}\footnote{It is worth noting that the \textit{global linear convergence} used in this paper, means that the sequence converges to a global optimum exponentially fast.} of FGD was developed in \cite{Sanghavi-FGD2016}.
Moreover, facing the challenge in large scale applications with a big $n$, stochastic algorithms \cite{Robbins-SGD1951} have been widely adopted nowadays, that is, at each iteration, we only use the gradient information of one or a small batch of the whole sample instead of the full gradient over $n$ samples. However, due to the variance of such stochastic gradients, the stochastic gradient descent (SGD) method only has a sublinear convergence rate even in the strongly convex case. Various variance reduction techniques have been proposed in literature (see, e.g., \cite{Johnson-Zhang-svrg2013,Bach-SAG}), which resume the linear convergence for strongly convex problems. However, it is still open whether these methods can be adapted to the nonconvex problem \eqref{Eq:SD-determ-ncvx} while enjoying the linear convergence to global optima.
%For problem \eqref{Eq:SD-stoch},
%the Stochastic Factored Gradient Descent (SFGD) algorithm can be naturally described as follows: at the $t$-th iteration, randomly pick an $i_t \in \{1,\ldots, n\}$, and then update the next iteration via
%\begin{align}
%\label{alg:SFGD}
%U^{t+1} = U^t - \eta_t \nabla f_{i_t}(X^t)U^t,
%\end{align}
%where
%$\eta_t >0$ is a diminishing step size.
%However, it is well known that the convergence rate of stochastic gradient descent (SGD) method is sublinear, even in the strongly convex case. This motivates us to study a faster stochastic algorithm for solving the nonlinear semidefinite optimization problem \eqref{Eq:SD-stoch}.

The main contribution in this paper is to fill in this gap by showing that, when adapted to the nonconvex problem \eqref{Eq:SD-determ-ncvx}, our proposed versions of stochastic variance reduced gradient (SVRG) method can find the global optimum of the original problem \eqref{Eq:SD-stoch} at a linear convergence rate when the initial choice lies in a prescribed neighbour of the global optimum and the objective function is restricted strongly convex. The initial choice condition here improves the one proposed for FGD in \cite{Sanghavi-FGD2016}. Moreover, our proposal includes both the fixed step sizes and the adaptive ones using a stabilized modification of Barzilai-Borwein (BB) step sizes \cite{BB-stepsize1988} which adapts to the non-convex problems when the curvature is not guaranteed as in strongly convex cases. Finally, experiments on both matrix sensing and ordinal embedding demonstrate the effectiveness of the proposed scheme. %To establish the global linear convergence rate, we have to introduce a completely different line of analysis in the vector setting, as well as modify existing proofs used in the deterministic setting.

%It is well known that the convergence rate of SVRG for the vector optimization (i.e., the argument is vector) is linear in the strongly convex case (see, \cite{Johnson-Zhang-svrg2013,Tan-Ma-bbsvrg2016}). However, its convergence is still unknown when applied to the semidefinite optimization problem \eqref{Eq:SD-stoch}. The main contribution of this paper is to fill in this gap.
%Under certain conditions, we establish the expected \textit{global linear convergence} of SVRG, starting from a good initial point. The critical role of initialization in the convergence guarantee is mainly due to SVRG is intrinsically designed for the \textit{nonconvex} problem \eqref{Eq:SD-determ-ncvx}. Thus, it is generally necessary to require a good initialization to find a global optimum.
%Several step-size strategies including the fixed, Barzilai-Borwein (BB) \cite{BB-stepsize1988}  and its variant called stabilized BB step sizes are considered and analyzed in a unified framework. A series of numerical experiments are also conducted to show the effectiveness of the proposed algorithm. To establish the global linear convergence rate, we have to introduce a completely different line of analysis in the vector setting, as well as modify existing proofs used in the deterministic setting.

%\subsection{Organization and Notations}

The reminder of this paper is organized as follows. Section \ref{sc:algorithms} introduce some algorithmic background with our proposal. Section \ref{sc:main result} presents the main convergence results.
%Section \ref{sc:relatedwork} provides the related works and some discussions.
Section \ref{sc:initial-scheme} provides some initialization schemes.
Section \ref{sc:proof} outlines the proof of our main theorem.
Section \ref{sc:experiment} provides some applications to verify our theoretical findings and show the effectiveness of the proposed algorithms.
We conclude this paper in Section \ref{sc:conclusion}.

{\bf Notations:} For any two matrices $X, Y \in \mathbb{R}^{p\times p}$, their inner product is defined as $\langle X, Y \rangle = \mathrm{tr}(X^TY)$. We denote $\mathbb{S}_+^p$ as the set of positive semidefinite matrices of size $p\times p$. For any matrix $X \in \mathbb{R}^{p\times p}$, $\|X\|_F$ and $\|X\|_2$ denote its Frobenius and spectral norms, respectively, and $\sigma_{\min}(X)$ and $\sigma_{\max}(X)$ denote the smallest and largest \textit{strictly positive} singular values of $X$,  denote $\tau(X):= \frac{\sigma_{\max}(X)}{\sigma_{\min}(X)}$, with a slight abuse of notation, we also use $\sigma_1(X)\equiv \sigma_{\max}(X) \equiv \|X\|_2$, and $X_r$ denotes the rank-$r$ approximation of $X$ via its truncated singular value decomposition (SVD) for any $r\leq p$. ${\bf I}_p$ denotes the identity matrix with the size $p\times p$. We will omit the subscript $p$ of ${\bf I}_p$ if there is no confusion in the context.

\section{Algorithms}
\label{sc:algorithms}
%We first introduce some algorithmic background for solving problem \eqref{Eq:SD-determ-ncvx}, and then present the global linear convergence of SVRG.

%\subsection{Algorithmic Background}

\begin{algorithm}[t]
{\small
\begin{algorithmic}\caption{SVRG for Problem \eqref{Eq:SD-stoch}}\label{alg:SVRG}
\STATE {\bf Parameters}: update frequency $m$, step size (or learning rate) $\{\eta_k\}$, initial point $\tilde{U}^0 \in \mathbb{R}^{p \times r}$
\smallskip
\FOR{k=0,1,\ldots}
\STATE $\tilde{X}^k := \tilde{U}^k {(\tilde{U}^k)}^T$
\STATE $g_k = \frac{1}{n} \sum_{i=1}^n \nabla f_i(\tilde{X}^k)\tilde{U}^k$
\STATE $U^0 = \tilde{U}^k$
\smallskip
\FOR{$t=0,\ldots,m-1$} %\textbf{do}:
\STATE $X^t = U^t {U^t}^T$
\STATE Randomly pick $i_t \in \{1,\ldots,n\}$
\STATE $U^{t+1} = U^t - \eta_k (\nabla f_{i_t}(X^t)U^t - \nabla f_{i_t}(\tilde{X}^k)\tilde{U}^k + g_k)$
\ENDFOR
\smallskip
\STATE $\tilde{U}^{k+1} = U^m$
\ENDFOR
\end{algorithmic}}
\end{algorithm}

Without loss of generality, we assume that $f$ is a symmetric function, i.e., $f(X)=f(X^T)$ throughout the paper. For $X=UU^T$, the gradient of $g(U):=f(UU^T)$ is
\[
\nabla g(U) = (\nabla f(UU^T)+\nabla f(UU^T)^T)U = 2 \nabla f(X)U.
\]

\textbf{A. FGD:} The FGD method proposed by \cite{Sanghavi-FGD2016} can be described as follows: let $U^t$ be the iterate at the $t$-th iteration and $X^t:=U^t {U^t}^T$, then the next iterate $U^{t+1}$ is updated according to the following
\begin{align}
\label{alg:FGD-fixed}
U^{t+1} = U^t - \eta \nabla f(X^t)U^t,
\end{align}
where $\eta >0$ is a step size.

\textbf{B. Stochastic FGD (SFGD):} As a stochastic counterpart of FGD \eqref{alg:FGD-fixed}, SFGD here can be naturally described as follows: at the $t$-th iteration, randomly pick an $i_t \in \{1,\ldots, n\}$, then update the next iteration via
\begin{align}
\label{alg:SFGD}
U^{t+1} = U^t - \eta_t \nabla f_{i_t}(X^t)U^t,
\end{align}
where
$\eta_t >0$ is a diminishing step size.

\textbf{C. SVRG:}\footnote{Besides SVRG, there are some other accelerated stochastic methods like SAG, SDCA and their variants. We focus on SVRG mainly due to SVRG might require less storage than SAG and SDCA, and thus it may be more suitable for the applications considered in this paper.}
The SVRG method was firstly proposed by \cite{Johnson-Zhang-svrg2013} for minimizing a finite sum of convex functions with a vector argument. The main idea of SVRG is adopting the variance reduction technique to accelerate SGD and achieve the linear convergence rate. Specifically, SVRG for solving problem \eqref{Eq:SD-determ-ncvx} can be described as in Algorithm \ref{alg:SVRG}. There are mainly two loops including an inner loop and an outer loop in SVRG.
%The inner loop adopted mainly aims to reduce the variance, and thus SVRG can overcome the main challenge of SGD as illustrated in \cite{Johnson-Zhang-svrg2013}.
One important implementation issue of SVRG is the tuning of the step size.
%In Algorithm \ref{alg:SVRG},
There are mainly two classes of step sizes: determined or data adaptive. Here we discuss three particular choices.
\begin{enumerate}
\item[(a)]
Fixed step size \cite{Johnson-Zhang-svrg2013}:
\begin{equation}
\label{Eq:fixed-stepsize}
\eta_k \equiv \eta, \quad \text{for some}\ \eta>0.
\end{equation}

\item[(b)]
Barzilai-Borwein (BB) step size \cite{BB-stepsize1988,Tan-Ma-bbsvrg2016}: given an initial $\eta_0 >0$ and for $k\geq 1$, let $\tilde{g}_k := \nabla f(\tilde{X}^k)$,
\begin{equation}
\label{Eq:bb-stepsize}
\eta_k = \frac{1}{m} \cdot \frac{\|\tilde{X}^k - \tilde{X}^{k-1}\|_F^2}{|\langle \tilde{X}^k - \tilde{X}^{k-1}, \tilde{g}_k - \tilde{g}_{k-1}\rangle|}.
\end{equation}

Note that such a BB step size is originally studied for strongly convex objective functions \cite{Tan-Ma-bbsvrg2016}, and it may be breakout if there is no guarantee of the curvature of $f$ like in nonconvex cases. In order to avoid such possible instability of \eqref{Eq:bb-stepsize} in our studies, a variant of BB step size, called the \textit{stabilized BB} step size, is suggested as follows.
\item[(c)]
Stabilized BB (SBB) step size: given an initial $\eta_0 >0$ and an $\epsilon \geq 0$, for $k \geq 1$,
\begin{align}
\label{Eq:sbb-stepsize}
&\eta_k = \frac{1}{m} \times\\
& \frac{\|\tilde{X}^k - \tilde{X}^{k-1}\|_F^2}{|\langle \tilde{X}^k - \tilde{X}^{k-1}, \tilde{g}_k - \tilde{g}_{k-1}\rangle| + \epsilon \|\tilde{X}^k - \tilde{X}^{k-1}\|_F^2}. \nonumber
\end{align}
%for some small $\epsilon>0.$
\end{enumerate}
%The original SVRG introduced in \cite{Johnson-Zhang-svrg2013} uses a fixed step size \eqref{Eq:fixed-stepsize}. We still name it as \textbf{SVRG} in this paper.
%The BB step size \eqref{Eq:bb-stepsize} was firstly incorporated into SVRG by \cite{Tan-Ma-bbsvrg2016} to avoid tuning the step size by hand. We call the SVRG method with BB step size \eqref{Eq:bb-stepsize} as \textbf{SVRG-BB} henceforth. However, as demonstrated by \cite{Tan-Ma-bbsvrg2016,Ma-SVRG-SBB2017}, the BB step size is efficient for the strongly convex function, yet it might be unstable when applied to the non-strongly convex function. This is mainly due to that the denominator in the BB step size \eqref{Eq:bb-stepsize} might be close to \textit{zero} or even \textit{negative} in the non-strongly convex (particularly, nonconvex) case. Thus, in \cite{Ma-SVRG-SBB2017}, the stabilized BB step size was introduced for SVRG to avoid the possible instability of BB step size when applied to practical applications like the ordinal embedding problem, of which the objective function might be non-strongly convex.
Throughout the rest of paper, with a slight abuse, we still name the original SVRG with a fixed step size as \textbf{SVRG}, and call the SVRG with stabilized BB step size \eqref{Eq:sbb-stepsize} as \textbf{SVRG-SBB}$_\epsilon$, and particularly, we call SVRG with BB step size as \textbf{SVRG-SBB}$_0$.
Besides the above three step sizes, there are some other schemes like the diminishing step size and the use of smoothing technique in BB step size as discussed in \cite{Tan-Ma-bbsvrg2016}. However, we mainly focus on the listed three step sizes in this paper due to they have been demonstrated to be effective in practice. Moreover, we only consider the Option-I suggested in \cite{Johnson-Zhang-svrg2013} for Algorithm \ref{alg:SVRG}, since Option-I in SVRG is generally a more natural and better choice than Option-II as demonstrated in both \cite{Johnson-Zhang-svrg2013} and \cite{Tan-Ma-bbsvrg2016} in the vector setting.

\section{Global Linear Convergence of SVRGs}
\label{sc:main result}

To present our main convergence results, we need the following assumptions.

\begin{Assump}
\label{Assump:objfun}
% Let $f(X) = \frac{1}{n}\sum_{i=1}^n f_i(X)$.
Each $f_i$ ($i=1,\ldots,n$)  satisfies the following:
\begin{enumerate}
\item[(a)]
$f_i$ is $L$-Lipschitz differentiable for some constant $L>0$, i.e., $f_i$ is smooth and $\nabla f_i$ is Lipschitz continuous satisfying
\[
\|\nabla f_i(X) - \nabla f_i(Y)\|_F \leq L \|X-Y\|_F, \ \forall X, Y \in \mathbb{S}_+^p.
\]

\item[(b)]
$f_i$ is $(\mu, r)$-restricted strongly convex for some constants $\mu>0$ and $r \leq p$, i.e., for any $X, Y \in \mathbb{S}_+^p$ with rank-$r$
\[
f_i(Y) \geq f_i(X) + \langle \nabla f_i(X), Y-X\rangle + \frac{\mu}{2} \|Y-X\|_F^2.
\]
\end{enumerate}
\end{Assump}

Assumption \ref{Assump:objfun} implies that $f$ is also $L$-Lipschitz differentiable and $(\mu, r)$-restricted strongly convex. For any $L$-Lipschitz differentiable and $(\mu, r)$-restricted strongly convex function $h$, the following hold (\cite{Nestrov-2004})
\begin{align*}
& h(Y) \leq h(X) + \langle \nabla h(X), Y-X \rangle + \frac{L}{2}\|Y-X\|_F^2, \\
%\label{Eq:Lip-ineq}\\
& \mu \|X-Y\|_F^2 \\
&\leq \langle \nabla h(X) - \nabla h(Y), X-Y \rangle \leq L \|X-Y\|_F^2,
%\label{Eq:innerproduc-ineq}
\end{align*}
where the first inequality holds for any $X, Y \in \mathbb{S}_+^p$, and the second inequality holds for any $X, Y \in \mathbb{S}_+^p$ with rank $r$, the first inequality and the right-hand side of the second inequality hold for the Lipschitz continuity of $\nabla h$, and the left-hand side of the second inequality is due to the $(\mu, r)$-restricted strong convexity of $h$.

Let $X^*$ be a global optimum of problem \eqref{Eq:SD-stoch} with rank $r^* : = \mathrm{rank}(X^*)$, $X_r^*$ be its rank-$r$ ($r\leq r^*$) best approximation via truncated singular value decomposition (SVD), and $U_r^*$ be a decomposition of $X_r^*$ via $X_r^* = U_r^* {U_r^*}^T$.
Under Assumption \ref{Assump:objfun}, we define the following constants:
\begin{align}
& \kappa := \frac{L}{\mu}, \quad \gamma_0 := \frac{2(\sqrt{2}-1)}{3\kappa}, \label{Eq:kappa}\\
%& \label{gamma0}\\
& \bar{\eta} := \min\left\{\frac{(1-\sqrt{\gamma_0})^2}{\frac{\|\nabla f(X_r^*)\|_F}{L \sigma_r(X_r^*)} +(2\sqrt{\gamma_0}+ \gamma_0)\tau(U_r^*)},1 \right\}, \label{bar eta}\\
& \xi := \bar{\eta}(1-\bar{\eta}/2), \label{xi}
\end{align}
where $\tau(X_r^*) : = \frac{\sigma_1(X_r^*)}{\sigma_r(X_r^*)}$ and
%$\sigma_r(X_r^*)$ is the $r$-th largest singular value of $X_r^*$,
$\tau(U_r^*) : = \frac{\sigma_1(U_r^*)}{\sigma_r(U_r^*)}$.
%is the condition number of $U_r^*$,
$\kappa \geq 1$ is generally called the \textit{condition number} of the objective function.  Thus, $0<\gamma_0 \leq \frac{2(\sqrt{2}-1)}{3}$
%and
%\begin{align*}
%&\min\left\{\frac{\left(1-\sqrt{\frac{2(\sqrt{2}-1)}{3}}\right)^2}{\left[\frac{\|\nabla f(X_r^*)\|_F}{L \cdot \|X_r^*\|} + \frac{2\sqrt{\frac{2(\sqrt{2}-1)}{3}}+ \frac{2(\sqrt{2}-1)}{3}}{\tau(U_r^*)}\right]\cdot \tau(X_r^*)},1 \right\}  \\
%&:= \tilde{\eta}\leq\bar{\eta} \leq 1,
%\end{align*}
and
$
0<\xi \leq 1/2.
$

As $r$ is used in the alternative nonconvex problem \eqref{Eq:SD-determ-ncvx}, the sequence $\{\tilde{X}^k\}$ generated by SVRG in Algorithm \ref{alg:SVRG} is at least rank-$r$, and can only converge to a rank-$r$ matrix if it is convergent. Therefore, we impose the following assumption to guarantee that the distance between $X_r^*$ and $X^*$ should be relatively small, otherwise, the introduced problem \eqref{Eq:SD-determ-ncvx} is not a good alterative of the original problem \eqref{Eq:SD-stoch}.

\begin{Assump}[rank-$r$ approximation error]
\label{Assump:r-approx error}
Let $X^*$ be a global optimum of problem \eqref{Eq:SD-stoch}, $X_r^*$ be the rank-$r$ approximation of $X^*$ for a given positive integer $r \leq r^*:=\mathrm{rank}(X^*)$. The following holds
\[
\|X_r^* - X^*\|_F < \frac{\sqrt{2}-1}{\sqrt{3}}\xi^{1/2} \kappa^{-1}\cdot\sigma_r(X^*),
\]
where $\kappa$ is specified in \eqref{Eq:kappa}, and $\sigma_r(X^*)$ is the $r$-th largest singular value of $X^*$.
\end{Assump}

Assumption \ref{Assump:r-approx error} is a regular assumption used in literature (say, \cite{Sanghavi-FGD2016}).
Roughly speaking, Assumption \ref{Assump:r-approx error} can be regarded as some noise assumption on problem \eqref{Eq:SD-stoch}. On the other hand, Assumption \ref{Assump:r-approx error} is imposed to guarantee the uniqueness of the rank-$r$ best approximation $X_r^*$. Otherwise, when $\|X_r^* - X^*\|_F \geq \sigma_r(X^*)$, if $X^* = {\bf I}_5$, i.e., an identity matrix with the size $5\times 5$, then $X_r^*$ with $r=4$ has five possible candidates.
Such assumption naturally holds for $r=r^*$, and when $r<r^*$, it might be satisfied if the singular values of $X^*$ possess certain \textit{compressible property}\footnote{$\sigma_i(X^*)$ decays in a power law, i.e., $\sigma_i(X^*) \leq C i^{-q}, i=1, 2, \ldots, p$ for some constants $C, q>0$}.
Under Assumption \ref{Assump:r-approx error}, we define several positive constants as follows: $\Delta : = \frac{(\sqrt{2}-1)^2 \xi^2\sigma_r^2(X_r^*)}{3\kappa^2} - \xi \|X_r^* - X^*\|_F^2, $
$\tilde{\Delta} : = \frac{4(\sqrt{2}-1)^2 \xi^2\sigma_r^2(X_r^*)}{9\kappa^2} - \xi \|X_r^* - X^*\|_F^2, $
\begin{align}
%&\Delta : = \frac{(\sqrt{2}-1)^2 \xi^2\sigma_r^2(X_r^*)}{3\kappa^2} - \xi \|X_r^* - X^*\|_F^2, \label{Delta}\\
%&\tilde{\Delta} : = \frac{4(\sqrt{2}-1)^2 \xi^2\sigma_r^2(X_r^*)}{9\kappa^2} - \xi \|X_r^* - X^*\|_F^2, \label{Delta1}\\
&\gamma_l :=\frac{2(\sqrt{2}-1)\xi\sigma_r(X_r^*)}{3\kappa} - \sqrt{\Delta}, \label{gamma l} \\
&\gamma_u :=\frac{2(\sqrt{2}-1)\xi\sigma_r(X_r^*)}{3\kappa} + \sqrt{\Delta}, \label{gamma u}\\
&\tilde{\gamma}_l :=\frac{2(\sqrt{2}-1)\xi\sigma_r(X_r^*)}{3\kappa} - \sqrt{\tilde{\Delta}}, \label{tilde_gamma l}\\
&\tilde{\gamma}_u :=\frac{2(\sqrt{2}-1)\xi\sigma_r(X_r^*)}{3\kappa} + \sqrt{\tilde{\Delta}}. \label{tilde_gamma u}
\end{align}
Note that the following relations hold
\begin{align}
&\gamma_l + \gamma_u = \frac{4(\sqrt{2}-1)\xi\sigma_r(X_r^*)}{3\kappa}, \label{Eq:gamma-l-u}\\
%&\gamma_l\cdot \gamma_u = \xi \|X_r^* - X^*\|_F^2 + \frac{(\sqrt{2}-1)^2 \xi^2\sigma_r^2(X_r^*)}{9\kappa^2},\label{Eq:gamma-l-u1}\\
%\end{align}
%and $\tilde{\gamma}_l$ and $\tilde{\gamma}_u$ satisfy the following relations:
%\begin{align}
%&\tilde{\gamma}_l + \tilde{\gamma}_u = \frac{4(\sqrt{2}-1)\xi\sigma_r(X_r^*)}{3\kappa}, \label{Eq:gamma1-l-u}\\
%&\tilde{\gamma}_l\cdot \tilde{\gamma}_u = \xi \|X_r^* - X^*\|_F^2, \label{Eq:gamma1-l-u1}\\
& \tilde{\gamma}_l < \gamma_l < \gamma_u < \tilde{\gamma}_u \leq \gamma_0 \sigma_r(X_r^*), \label{Eq:gamma2-l-u}
%& \gamma_u < \gamma_0 \sigma_r(X_r^*). \label{Eq:gamma-u-0}
\end{align}
where the last inequality of \eqref{Eq:gamma2-l-u} holds for $0<\xi \leq 1/2$ and $\tilde{\gamma}_u \leq 2\xi \gamma_0 \sigma_r(X_r^*) \leq \gamma_0 \sigma_r(X_r^*)$.

We also need the following common assumption on the stochastic direction, which has been widely used in literature on stochastic algorithms (say, \cite{Bottou-SGD2016} and reference therein).
\begin{Assump}[Unbiasedness]
\label{Assump:moment}
 $\{\nabla f_{i_t}(X^t)U^t\}$ satisfies $\mathbb{E}_{i_t}[\nabla f_{i_t}(X^t)U^t] = \nabla f(X^t)U^t, \ \forall t\in \mathbb{N}$.
\end{Assump}

If $i_t$ is uniformly sampled,
(see \cite{Needell-SGD2014,Zhang-SGDImportSamp2015} for studies on importance sampling),
then the above assumption can be satisfied.
Under Assumptions \ref{Assump:objfun}-\ref{Assump:moment}, let ${\cal N}_{\gamma_0}:= \{U:\|U-U_r^*\|^2_F \leq \gamma_0 \sigma_r(X_r^*)\}$, and we define the following constants: ${\cal B}:= \sup_{U \in {\cal N}_{\gamma_0}} \|UU^T\|_F, $
$B_0 := \sup_{U \in {\cal N}_{\gamma_0}} \left\{\mathbb{E}_{i_t}[\|\nabla f_{i_t}(UU^T)\|_F^2] - \|\nabla f(UU^T)\|_F^2\right\}, $
$B_1 := \sup_{U \in {\cal N}_{\gamma_0}} \|\nabla f(UU^T)\|^2_F, $
\begin{align}
%&{\cal B}:= \sup_{U \in {\cal N}_{\gamma_0}} \|UU^T\|_F, \label{Eq:B-cal}\\
%&B_0 := \nonumber\\
%&\sup_{U \in {\cal N}_{\gamma_0}}
%\left\{\mathbb{E}_{i_t}[\|\nabla f_{i_t}(UU^T)\|_F^2] - \|\nabla f(UU^T)\|_F^2\right\}, \label{Eq:B0}\\
%&B_1 := \sup_{U \in {\cal N}_{\gamma_0}} \|\nabla f(UU^T)\|^2_F, \label{Eq:B1}\\
& B_2 := 4\left[ 2L^2 {\cal B} ({\cal B} + \|X_r^*\|_F) + B_0 + B_1\right], \label{Eq:B2}\\
%&\delta : = \sqrt{\tilde{\Delta}}+\sqrt{\Delta}, \label{Eq:delta}\\
& \theta: = \frac{2\xi B_2}{L(\sqrt{\tilde{\Delta}}- \sqrt{\Delta})} = \frac{18B_2 \kappa \delta}{(\sqrt{2}-1)^2 \xi \mu \sigma_r^2(X_r^*)}, \label{Eq:theta}\\
& \eta_{\max} := \min \left\{\zeta_1, \zeta_2,\frac{1}{2\theta}\right\}, \label{Eq:eta-max}
\end{align}
where $\delta : = \sqrt{\tilde{\Delta}}+\sqrt{\Delta}, $ $\zeta_1 := \frac{1}{12\left[2L \cdot \kappa {\cal B} + \frac{B_0 + B_1}{(\sqrt{2}-1)\mu \sigma_r(X_r^*)}\right]},$
and $\zeta_2 := \frac{(\sqrt{2}-1)\mu \xi \sigma_r(X_r^*)}{12B_2}$. It can be seen that ${\cal B}$ is the upper bound of $X = UU^T$, $B_0$ represents variance of the stochastic gradient of $f$, and $B_1$ is the upper bound of the squared Frobenius norm of gradient $\nabla f(UU^T)$, restricted to the closed ball ${\cal N}_{\gamma_0}$.
%From the above definitions, there hold
%\begin{align*}
%& {\cal B} \leq \sup_{U \in {\cal N}_{\gamma_0}} \|UU^T - U_r^* {U_r^*}^T\|_F + \|X_r^*\|_F \\
%& \mathop{\leq}^{\text{(Lemma \ref{Lemm:X-Xr*})}} (2\sqrt{\gamma_0}+\gamma_0)\cdot\tau(U_r^*)\cdot \sigma_r(X_r^*) + \|X_r^*\|_F,\\
%& B_1 = \sup_{U \in {\cal N}_{\gamma_0}} \|\nabla f(UU^T) \|_F^2 \\
%& \leq 2L^2({\cal B}+\|X_r^*\|_F)^2+2\|\nabla f(X_r^*)\|_F^2.
%\end{align*}

%\subsection{Convergence Results}

Let $\{\eta_k\}$ be a sequence satisfying $\eta_k \in (0,\eta_{\max})$ for any $k \in \mathbb{N}$. Given a positive integer $m$, define
\begin{align}
%1- \frac{\eta_k L (\sqrt{\tilde{\Delta}}-\sqrt{\Delta})}{2\xi} =
&\rho_k := 1-\frac{\eta_k (\sqrt{2}-1)^2 \xi \mu \sigma_r^2(X_r^*)}{18 \kappa \delta}, \label{Eq:rho}\\
&\tilde{\rho}_k: = \rho_k^m + (1-\rho_k^m) \eta_k \theta. \label{Eq:rho-tilde}
\end{align}
It is easy to check that $0<\rho_k<1$ and $0<\tilde{\rho}_k<1$.
Based on the above defined constants, we present our main theorem as follows.

\begin{Thm}[Linear convergence of SVRG]
\label{Thm:svrg}
Let $\{\tilde{U}^k\}$ be a sequence generated by Algorithm \ref{alg:SVRG}. Suppose that Assumptions \ref{Assump:objfun}-\ref{Assump:moment} hold, and that $\eta_k \in (0,\eta_{\max})$. The following hold:
(a) if $\gamma_l < \|\tilde{U}^0-U_r^*\|_F^2<\gamma_u,$
%\begin{align}
%\label{Assump:initial-SFGD}
%\gamma_l < \|\tilde{U}^0-U_r^*\|_F^2<\gamma_u,
%\end{align}
there hold
\begin{enumerate}
\item[(a1)] $\{\mathbb{E}[\|\tilde{U}^{k} - U_r^*\|_F^2]\}$ is monotonically decreasing,

\item[(a2)](\textbf{Linear convergence}) for any $k\geq 1,$
\begin{align}
&\mathbb{E}[\|\tilde{U}^{k} - U_r^*\|_F^2]\label{Eq:linearconv-svrg}\\
&\leq \left( \prod_{i=0}^{k-1}\tilde{\rho}_i \right) \cdot \|\tilde{U}^0 - U_r^*\|_F^2 + \tilde{\gamma_l} \times\nonumber\\
&  \left[ \sum_{t=0}^{k-2} \left(\prod_{i=t+1}^{k-1} \tilde{\rho_i} \cdot (1-(\rho_t)^m)\right) + \left(1- (\rho_{k-1})^m \right) \right]. \nonumber
\end{align}
\end{enumerate}
(b) In addition, if $\|\tilde{U}^0 - U_r^*\|^2_F \leq \gamma_l$, then $\mathbb{E}[\|\tilde{U}^{k}- U_r^*\|_F^2] \leq \gamma_l$ for any $k\in \mathbb{N}$.
\end{Thm}

\begin{figure}[!t]
\begin{minipage}[b]{0.98\linewidth}
\centering
\includegraphics*[scale=0.3]{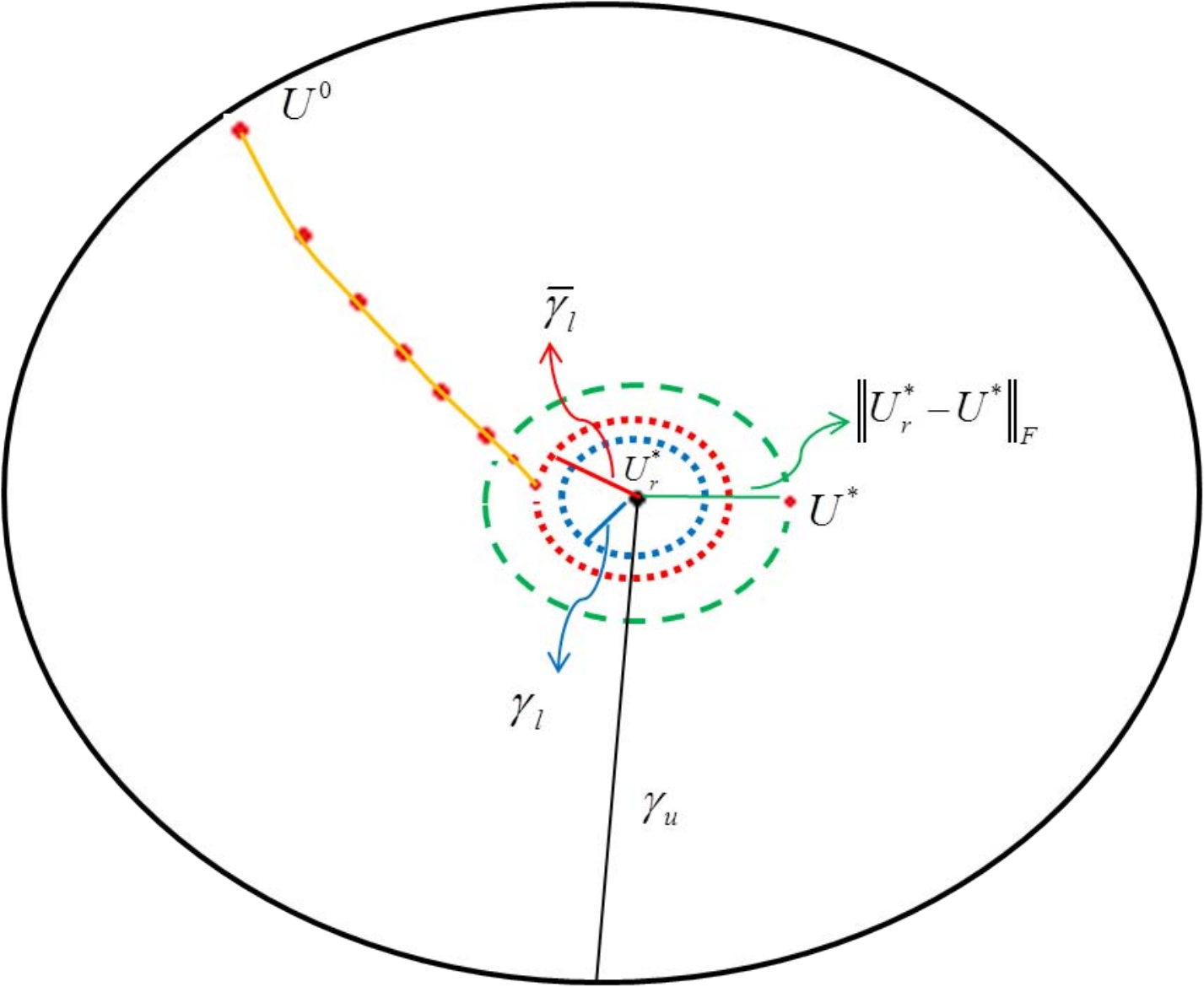}
\end{minipage}
\caption{ Convergence path of SVRG.
}
\label{Fig:convergence_path_svrg}
\end{figure}

The above theorem holds for a generic step size satisfying $\eta_k \in (0,\eta_{\max})$.
Actually, if $\{\eta_k\}$ is lower bounded by a positive constant $\eta_{\min}$ and obviously,
%there is a lower bound of $\eta_k$, i.e., there exists a positive
$\eta_{\min} < \eta_{\max}$,
%such that $\eta_k \in (\eta_{\min},\eta_{\max})$ for any $k\in \mathbb{N}$,
then by \eqref{Eq:rho} and \eqref{Eq:rho-tilde}, $\rho_k \in (\rho_{\min}, \rho_{\max})$ and
$\tilde{\rho}_k \in (\theta\eta_{\min}, \tilde{\rho}_{\max}),$
where $\rho_{\min}:=1-\frac{\eta_{\max}(\sqrt{2}-1)^2  \xi \mu \sigma_r^2(X_r^*)}{18 \kappa \delta}$, $\rho_{\max}:=1-\frac{\eta_{\min}(\sqrt{2}-1)^2  \xi \mu \sigma_r^2(X_r^*)}{18 \kappa \delta}$, and $\tilde{\rho}_{\max}:= \rho_{\max}^m + (1-\rho_{\max}^m)\theta\eta_{\max}<1$.
Thus, $\prod_{i=0}^{k-1}\tilde{\rho}_i \leq \left(\tilde{\rho}_{\max}\right)^k,$ and
\begin{align*}
& \sum_{t=0}^{k-2} \left(\prod_{i=t+1}^{k-1} \tilde{\rho_i} \cdot (1-(\rho_t)^m)\right) + (1- (\rho_{k-1})^m) \\
& \leq (1-\rho_{\min}^m)\cdot \left[ 1+ \sum_{t=0}^{k-2} (\tilde{\rho}_{\max})^{k-t-1}\right]\\
& = (1-\rho_{\min}^m)\cdot \frac{1-(\tilde{\rho}_{\max})^k}{1-\tilde{\rho}_{\max}}
= \frac{1-\rho_{\min}^m}{1-\rho_{\max}^m} \cdot \frac{1-(\tilde{\rho}_{\max})^k}{1-\theta \eta_{\max}}.
%: = c_{\max} \left[ 1-(\tilde{\rho}_{\max})^k \right].
\end{align*}
Let
$\bar{\gamma}_l: = \frac{1-\rho_{\min}^m}{1-\rho_{\max}^m} \cdot \frac{\tilde{\gamma}_l}{1-\theta \eta_{\max}}$.
According to the above two inequalities, \eqref{Eq:linearconv-svrg} implies that
\begin{align*}
\mathbb{E}[\|\tilde{U}^{k} - U_r^*\|_F^2] - \bar{\gamma}_l \leq \left(\tilde{\rho}_{\max}\right)^k (\|\tilde{U}^0 - U_r^*\|_F^2 - \bar{\gamma}_l),
\end{align*}
which shows the linear convergence of SVRG.
%As a consequence, Theorem \ref{Thm:svrg} shows some linear global convergence of SVRG.
Thus, Theorem \ref{Thm:svrg} shows certain \textit{global linear} convergence of SVRG, that is, the convergence to a global optimum starting from some good initial point, as depicted in Figure \ref{Fig:convergence_path_svrg}.
From Figure \ref{Fig:convergence_path_svrg},
starting from an initialization lying in a $\gamma_u$-neighborhood of $U_r^*$, SVRG converges exponentially fast until achieving a small $\bar{\gamma}_l$-neighborhood of $U_r^*$;
%, then stagnates and never jumps out from this smaller neighborhood.
%More specifically, Theorem \ref{Thm:svrg} demonstrates that when the initialization lies in the band between a larger $\gamma_u$-ball and a smaller $\gamma_l$-ball of $U_r^*$, then SVRG will converge exponentially fast to a $\bar{\gamma}_l$-ball of $U_r^*$ in expectation,
while if the initialization lies in the $\gamma_l$-ball of $U_r^*$, then SVRG will never escape from this small ball in expectation.

%However, it is generally hard to compare directly the values of $\gamma_l$ and $\bar{\gamma}_l$ because both depend on several terms. Moreover, it is worth noting that $\bar{\gamma}_l$ is only an upper bound of the radius of the ``limiting'' ball.

The comparisons on convergence results between FGD \cite{Sanghavi-FGD2016} and SVRG in the restricted strongly convex case are shown in Table \ref{Tab:comp_SVRG}. The convergence result of SVRG is presented in expectation. From Table \ref{Tab:comp_SVRG},
the requirement on the rank-$r$ approximation error can be relaxed from the order ${\cal O}(\frac{\sigma_r(X_r^*)}{\kappa^{1.5}\tau(X_r^*)})$ to ${\cal O}(\frac{\sigma_r(X_r^*)}{\kappa})$, and the requirement on the radius of initialization can be relaxed from ${\cal O}(\frac{\sigma_r(X_r^*)}{\kappa^{2}\tau^2(X_r^*)})$ to ${\cal O}(\frac{\sigma_r(X_r^*)}{\kappa})$, where $\kappa$ is the ``condition number'' of the objective function $f$ (specified in \eqref{Eq:kappa}), $\sigma_r(X_r^*)$ and $\tau(X_r^*)$ are respectively the $r$-th largest singular value and the condition number of the rank-$r$ approximation $X_r^*$ of the optimum $X^*$ with $r\leq r^* := \mathrm{rank}(X^*)$.

%the conditions on rank-$r$ approximation error and initialization of SVRG are comparable with those of SFGD developed in \cite{Zeng-SFGD2017} but much weaker than those of FGD established in \cite{Sanghavi-FGD2016}, in terms of the order. Moreover, the rate of SVRG is linear but that of SFGD is sublinear, i.e., ${\cal O}(1/k)$. Although FGD has the similar linear rate, we will show in the latter numerical experiments that SVRG converges much faster than FGD in terms of \textit{epoch number} (explained latter).

{
    \renewcommand\baselinestretch{.9}
    \renewcommand\arraystretch{0.6}
\begin{table}[ht!]
\tabulinesep=1.5mm
\centering
\small
%\resizebox{0.9\textwidth}{!}
\begin{tabu}{|l|c|c|}\hline
\tiny Algorithm      & \tiny $\|X^* - X_r^*\|_F$  & \tiny Initialization          \\ \hline
\tiny FGD (\cite{Sanghavi-FGD2016}) & \tiny ${\cal O}\left(\frac{\sigma_r(X_r^*)}{\kappa^{1.5}\tau(X_r^*)}\right)$ & \tiny ${\cal O}\left(\frac{\sigma_r(X_r^*)}{\kappa^2 \tau^2(X_r^*)}\right)$ \\ \hline
%\tiny SFGD (\cite{Zeng-SFGD2017})  & \tiny $\frac{\sigma_r(X_r^*)}{\kappa \tau^{0.5}(X_r^*)}$  & \tiny $\frac{\sigma_r(X_r^*)}{\kappa \tau(X_r^*)}$  & \tiny $\frac{\kappa \|X_r^* - X^*\|_F^2}{\sigma_r(X_r^*)}$  & \tiny $1/k$\\ \hline
\tiny SVRG (our)  & \tiny $\textcolor{red}{{\cal O}\left(\frac{\sigma_r(X_r^*)}{\kappa \tau^{0.5}(X_r^*)}\right)}$  & \tiny $\textcolor{red}{{\cal O}\left(\frac{\sigma_r(X_r^*)}{\kappa \tau(X_r^*)}\right)}$  \\ \hline
\end{tabu}\\
\caption{\small{Comparisons on convergence results (in order) between FGD \cite{Sanghavi-FGD2016} and SVRG (this paper) in the restricted strongly convex case.
%The convergence results of both SFGD and SVRG are shown in the sense of expectation.
%The order of the convergence conditions and results are presented in this table.
%We omit the notation ${\cal O}$ for the sake of limited space.
%The second and third columns show the requirements on the rank-$r$ approximation error and initialization.
%The fourth column shows the limiting recovery errors, which are measured by $\mathop{\lim}_{k\rightarrow \infty} \|U^k- U_r^*\|_F^2$.
%The last column shows the convergence rates, where $\hat{\rho}$ and $\tilde{\rho}_{\max}$ are positive constants less than 1.
}
}
\label{Tab:comp_SVRG}
\end{table}
}

In the following, we give a corollary to show the convergence of SVRG when adopting the considered three step-size strategies \eqref{Eq:fixed-stepsize}-\eqref{Eq:sbb-stepsize}.

\begin{Coro}[Convergence for different step sizes]
\label{Coro:3stepsizes}
Under conditions of Theorem \ref{Thm:svrg}, all claims in Theorem \ref{Thm:svrg} hold, if one of the following conditions holds:
\begin{enumerate}
\item[(1)]
$\eta \in (0,\eta_{\max})$ when a fixed step size is adopted;

%\item[(2)]
%$m > \frac{1}{\mu \eta_{\max}}$ when BB step size is adopted; \commyy{remove it}

\item[(2)]
$m > \frac{1}{(\mu+\epsilon) \eta_{\max}}$ for any $\epsilon\geq 0$ when SBB step size is adopted.
\end{enumerate}
\end{Coro}

%\commyy{we just need to discuss $SBB_\epsilon$.}
By the definition of SBB step size \eqref{Eq:sbb-stepsize}, and Assumption \ref{Assump:objfun}, we have
\[
\frac{1}{m(L+\epsilon)} \leq \eta_k \leq \frac{1}{m(\mu + \epsilon)}.
\]
%Thus, if $m >  \frac{1}{\mu \eta_{\max}}$, then BB step size $\eta_k < \eta_{\max}$ for any $k \in \mathbb{N}$.
%%Moreover, it is easy to show that $\eta_k$ is uniformly lower bounded by $\frac{1}{mL}$.
%Similarly, we can also show that the stabilized BB step size $\eta_k$ satisfies
%\[
%\frac{1}{m(L+\epsilon)} \leq \eta_k \leq \frac{1}{m(\mu + \epsilon)},
%\]
%and
Thus, if $m > \frac{1}{(\mu + \epsilon) \eta_{\max}}$, then $\eta_k < \eta_{\max}$ for any $k \in \mathbb{N}$.

From \eqref{gamma l}-\eqref{tilde_gamma l}, if $r=r^*$ then $\|X_r^* - X^*\|_F =0$, and thus $\tilde{\gamma}_l =0$ and $\gamma_u = \frac{(2+\sqrt{3})\cdot(\sqrt{2}-1)\xi\sigma_r(X_r^*)}{3\kappa}$.
However, in this case, $\gamma_l = \frac{(2-\sqrt{3})\cdot(\sqrt{2}-1)\xi\sigma_r(X_r^*)}{3\kappa}>0$. Thus, we cannot claim the exact recovery of a global optimum directly from Theorem \ref{Thm:svrg} even if $\|\tilde{U}^0 - U_r^*\|^2_F \leq \gamma_l$.
To circumvent this problem, we use a more consecutive step size, and get the following corollary showing the exact recovery of SVRG.
%we have the following corollary showing the exact recovery of SVRG but with a more
Let
\begin{align}
\label{Eq:bar-eta-max}
\bar{\eta}_{\max}:= \min\left\{\frac{L\gamma_u}{2B_2 \xi}, \eta_{\max}\right\}.
\end{align}

\begin{Coro}[Exact recovery when $r=r^*$]
\label{Coro:exact-recovery}
Let $\{\tilde{U}^k\}$ be a sequence generated by Algorithm \ref{alg:SVRG}. Let Assumptions \ref{Assump:objfun} and \ref{Assump:moment} hold. If the following conditions hold: (a) $r=r^*$, (b) $\eta_k \in (0,\bar{\eta}_{\max})$, and (c) $\|\tilde{U}^0-U_r^*\|_F^2<\frac{(2+\sqrt{3})\cdot(\sqrt{2}-1)\xi\sigma_r(X_r^*)}{3\kappa}$, then SVRG exactly recover the global optimum $X^*$ in expectation at a linear rate.
\end{Coro}
Corollary \ref{Coro:exact-recovery} shows that if fortunately, we can take $r$ as the exact rank $r^*$ of the global optimum, then SVRG can exactly find the global optimum in expectation exponentially fast, as long as the initialization lies in a neighborhood of the global optimum.
The proof of this corollary is presented in (Supplementary Material: Section 2.2).

\section{On Initialization Schemes}
\label{sc:initial-scheme}
According to our main theorem (see, Theorem \ref{Thm:svrg}),
the initialization should be close to $U_r^*$ to get the provable convergence.
%However, note that the studied algorithms are designed for solving a \textit{nonconvex} problem \eqref{Eq:SD-determ-ncvx}. This brings some challenges to get a good initialization.
In the following, we discuss some potential initialization schemes.

\textbf{Scheme I:} One common way is to use one of the standard convex algorithms (say, projected gradient descent method) and obtain a good initialization $U^0$, then switch to SVRG to get a higher precision solution. A specific implementation of this idea has been used in \cite{Recht-2016} to deal with the matrix sensing problem, and some theoretical guarantees of this scheme have been developed in \cite{Sanghavi-FGD2016}.
%However, this scheme is generally time-consuming and not scalable.

\textbf{Scheme II:}
%We describe another initialization scheme constructed from $0$, which has been used in many recent work (say, \cite{Sanghavi-FGD2016,Sanghavi2013}).
%More specifically, we first
Another way is firstly to get $X^0:=\frac{1}{\|\nabla f(0) - \nabla f({\bf e}_1 {\bf e}_1^T)\|_F}\cdot \mathrm{Proj}_{\mathbb{S}_+^p}(-\nabla f(0))$, then take $U^0 \in \mathbb{R}^{p\times r}$ such that $U^0 {U^0}^T = X_r^0 $, where $X_r^0$ is the rank-$r$ best approximation of $X^0$ via SVD, and ${\bf e}_1 \in \mathbb{R}^p$ is the vector with 1 as the first component and 0 as the other components.
The effectiveness of such scheme is guaranteed by \cite[Corollary 12]{Sanghavi-FGD2016} when the objective function is \textit{well-conditioned}, i.e., has a small $\kappa$.
%, while when $\kappa$ is large, it might not work.

\textbf{Scheme III:} Note that the previous two schemes need at least one SVD, which might be prohibitive in large scale applications. To avoid such an issue, random initialization can be exploited which actually works well in many applications.

%\section{Proof of Theorem \ref{Thm:svrg}}
\section{Outline of Proofs}
\label{sc:proof}

To prove Theorem \ref{Thm:svrg}, we need the following key lemma, which gives an error estimate of the inner loop.
%For any matrix $U \in \mathbb{R}^{p\times r}$, let $Q_U$ be a basis of the column space of $U$. Denote ${\cal P}_U := Q_UQ_U^T$. Then ${\cal P}_U \cdot U = U$. For any matrix $Y \in \mathbb{R}^{p\times p}$, ${\cal P}_UY$ is a projection of $Y$ onto the subspace spanned by $X:= UU^T$.

\begin{Lemm}[A key lemma]
\label{Lemm:innerloop}
Let $\{U^t\}_{t=0}^{m}$ be the sequence at the $k$-th inner loop. Let Assumptions \ref{Assump:objfun}, \ref{Assump:r-approx error} and \ref{Assump:moment} hold. Let $\eta_k \in (0,\eta_{\max})$. If $\gamma_l<\mathbb{E}[\|\tilde{U}^k - U_r^*\|_F^2] < \gamma_u$, then the sequence
$\{\mathbb{E}[\|U^t - U_r^*\|_F^2]\}$ is monotonically decreasing for $t=0,\ldots,m$, and
\begin{align}
\label{Eq:recur-innerloop}
&\mathbb{E}_{i_t}[\|U^{t+1} - U_r^*\|_F^2] \leq \frac{\eta_k L}{2} \|X^* - X_r^*\|_F^2\\
&+ \|U^t - U_r^*\|_F^2 - \frac{2(\sqrt{2}-1)}{3}\eta_k \mu \sigma_r(X_r^*) \|U^t - U_r^*\|_F^2 \nonumber\\
& + \frac{\eta_k L}{2\xi} \|U^t - U_r^*\|_F^4 + \eta_k^2 B_2 \cdot \|\tilde{U}^k - U_r^*\|_F^2 \nonumber
\end{align}
where $B_2$ is specified in \eqref{Eq:B2};
while if $\mathbb{E}[\|\tilde{U}^k - U_r^*\|_F^2]\leq \gamma_l$, then $\mathbb{E}[\|U^t - U_r^*\|_F^2]\leq \gamma_l$ for any $t=0,\ldots,m.$
\end{Lemm}
{\bf The sketch proof of Lemma \ref{Lemm:innerloop}:}
We prove this lemma by induction. Specifically, we first show that if $\gamma_l < \mathbb{E}[\|U^t - U_r^*\|_F^2] < \gamma_u$, then $\mathbb{E}[\|U^{t+1} - U_r^*\|_F^2]\leq \mathbb{E}[\|U^t - U_r^*\|_F^2]  < \gamma_u$ for $t=0,\ldots,m-1$. Furthermore, $\mathbb{E}[\|U^{t+1} - U_r^*\|_F^2]$ can be estimated via noting that
\begin{align*}
&\mathbb{E}_{i_t}[\|U^{t+1} - U_r^*\|_F^2] \\
&= \|U^t - U_r^*\|_F^2 + \eta_k^2 \mathbb{E}_{i_t}[\|v_k^t\|_F^2] \\
&- 2\eta_k \mathbb{E}_{i_t}[\langle v_k^t, U^t - U_r^* \rangle],
\end{align*}
where $v_k^t = \nabla f_{i_t}(X^t)U^t - \nabla f_{i_t}(\tilde{X}^k)\tilde{U}^k + \nabla f(\tilde{X}^k)\tilde{U}^k$,
and then establish the bounds of both $\mathbb{E}_{i_t}[\|v_k^t\|_F^2]$ and $\mathbb{E}_{i_t}[\langle v_k^t, U^t - U_r^* \rangle]$ via two lemmas shown in (Supplementary Material: Lemma 2 and Lemma 3), respectively.
The specific proof of this lemma is presented in (Supplementary Material: Section A).

Based on Lemma \ref{Lemm:innerloop}, we show the proof of Theorem \ref{Thm:svrg} as follows.

%\begin{proof}[Proof of Theorem \ref{Thm:svrg}]
{\textbf{Proof of Theorem \ref{Thm:svrg}}:}
By Lemma \ref{Lemm:innerloop}, if $\gamma_l<\|\tilde{U}^0 - U_r^*\|_F^2 < \gamma_u$, then for any $k \in \mathbb{N}$ and $t=0,\ldots,m-1$,
\[
\mathbb{E}[\|\tilde{U}^{k+1} - U_r^*\|_F^2] \leq \mathbb{E}[\|U^t - U_r^*\|_F^2] \leq \mathbb{E}[\|\tilde{U}^{k} - U_r^*\|_F^2].
\]
From \eqref{Eq:recur-innerloop} and by the definitions of $\tilde{\gamma}_l$ and $\tilde{\gamma}_u$, at the $k$-th inner loop, there holds
\begin{align*}
& \mathbb{E}[\|U^{t+1} - U_r^*\|_F^2] - \tilde{\gamma}_l \\
& \leq \left[ 1 - \frac{\eta_k L}{2\xi}(\tilde{\gamma}_u - \mathbb{E}[\|U^t - U_r^*\|_F^2])\right]\times \\
& \left(\mathbb{E}[\|U^t - U_r^*\|_F^2] - \tilde{\gamma}_l \right) + \eta_k^2 B_2 \cdot \|\tilde{U}^k - U_r^*\|_F^2\\
& \leq \left[ 1 - \frac{\eta_k L}{2\xi}\left(\sqrt{\tilde{\Delta}}-\sqrt{\Delta}\right)\right]\cdot (\mathbb{E}[\|U^t - U_r^*\|_F^2] - \tilde{\gamma}_l) \\
&+ \eta_k^2 B_2 \cdot \|\tilde{U}^k - U_r^*\|_F^2\\
&:= \rho_k (\mathbb{E}[\|U^t - U_r^*\|_F^2] - \tilde{\gamma}_l) + \eta_k^2 B_2 \cdot \|\tilde{U}^k - U_r^*\|_F^2,
\end{align*}
where the second inequality holds for $\mathbb{E}[\|U^t - U_r^*\|_F^2] < \gamma_u$ and
$\eta_k < \eta_{\max} \leq \frac{(\sqrt{2}-1)\mu \xi \sigma_r(x_r^*)}{12 B_2}
\leq \frac{6}{(\sqrt{2}-1)\mu \sigma_r(X_r^*)} \leq \frac{2\xi}{L(\sqrt{\tilde{\Delta}} - \sqrt{\Delta})}.
$
By the above inequality, we have
\begin{align*}
&\mathbb{E}[\|\tilde{U}^{k+1} - U_r^*\|_F^2] -\tilde{\gamma}_l
= \mathbb{E}[\|U^m - U_r^*\|_F^2] -\tilde{\gamma}_l \\
&\leq (\rho_k)^m (\mathbb{E}[\|\tilde{U}^{k} - U_r^*\|_F^2] -\tilde{\gamma}_l) \\
& + \eta_k^2 B_2 \cdot \frac{1-(\rho_k)^m}{1-\rho_k}\cdot \mathbb{E}[\|\tilde{U}^k - U_r^*\|_F^2]\\
&\leq (\rho_k)^m (\mathbb{E}[\|\tilde{U}^{k} - U_r^*\|_F^2] -\tilde{\gamma}_l) \\
& + \eta_k\theta (1-(\rho_k)^m) \mathbb{E}[\|\tilde{U}^{k} - U_r^*\|_F^2],
\end{align*}
where the final inequality is due to the definition of $\rho_k$ \eqref{Eq:rho}, i.e., $\rho_k = 1- \frac{\eta_k B_2}{\theta}$ and $\theta$ is specified in \eqref{Eq:theta}.
Therefore,
\begin{align*}
&\mathbb{E}[\|\tilde{U}^{k+1} - U_r^*\|_F^2] \leq  (1-(\rho_k)^m)\tilde{\gamma}_l\\
& + \left[(\rho_k)^m + \eta_k\theta (1-(\rho_k)^m)\right] \cdot \mathbb{E}[\|\tilde{U}^{k} - U_r^*\|_F^2] \\
&:= \tilde{\rho}_k\mathbb{E}[\|\tilde{U}^{k} - U_r^*\|_F^2]+(1-(\rho_k)^m)\tilde{\gamma}_l.
\end{align*}
Based on the above inequality, we get \eqref{Eq:linearconv-svrg} via a recursive way, and thus complete the proof of this theorem.
%As a consequence,
%\begin{align*}
%&\mathbb{E}[\|\tilde{U}^{k} - U_r^*\|_F^2]
%\leq \left( \prod_{i=0}^{k-1}\tilde{\rho}_i \right) \cdot \|\tilde{U}^0 - U_r^*\|_F^2 \\
%&+ \left[ \sum_{t=0}^{k-2} \left(\prod_{i=t+1}^{k-1} \tilde{\rho_i} \cdot (1-(\rho_t)^m)\right) + (1- (\rho_{k-1})^m) \right]\cdot \tilde{\gamma_l}.
%\end{align*}
%Thus, we complete the proof of this theorem.
$\Box$
%\end{proof}

\section{Experiments}
\label{sc:experiment}
In this section, we present two application examples to show the effectiveness of the proposed algorithm and also verify our developed theoretical results.
%More experiment results of the ordinal embedding application can be found in (Supplementary Material: Section 2).

\subsection{Matrix Sensing}

We consider the following matrix sensing problem
\begin{align*}
\mathop{\mathrm{min}}_{X \succeq 0} \quad f(X) = \frac{1}{2n} \sum_{i=1}^n (b_i - \langle A_i,X \rangle)^2,
\end{align*}
where $X \in \mathbb{R}^{p\times p}$ is a low-rank matrix, $A_i \in \mathbb{R}^{p\times p}$ is a sub-Gaussian independent measurement matrix of the $i$-th sample, $b_i \in \mathbb{R}$, and $n \in \mathbb{N}$ is the sample size.

\begin{figure}[h]
%\vspace{-.3in}
%\centerline{\fbox{This figure intentionally left non-blank}}
\includegraphics[width=.47\linewidth]{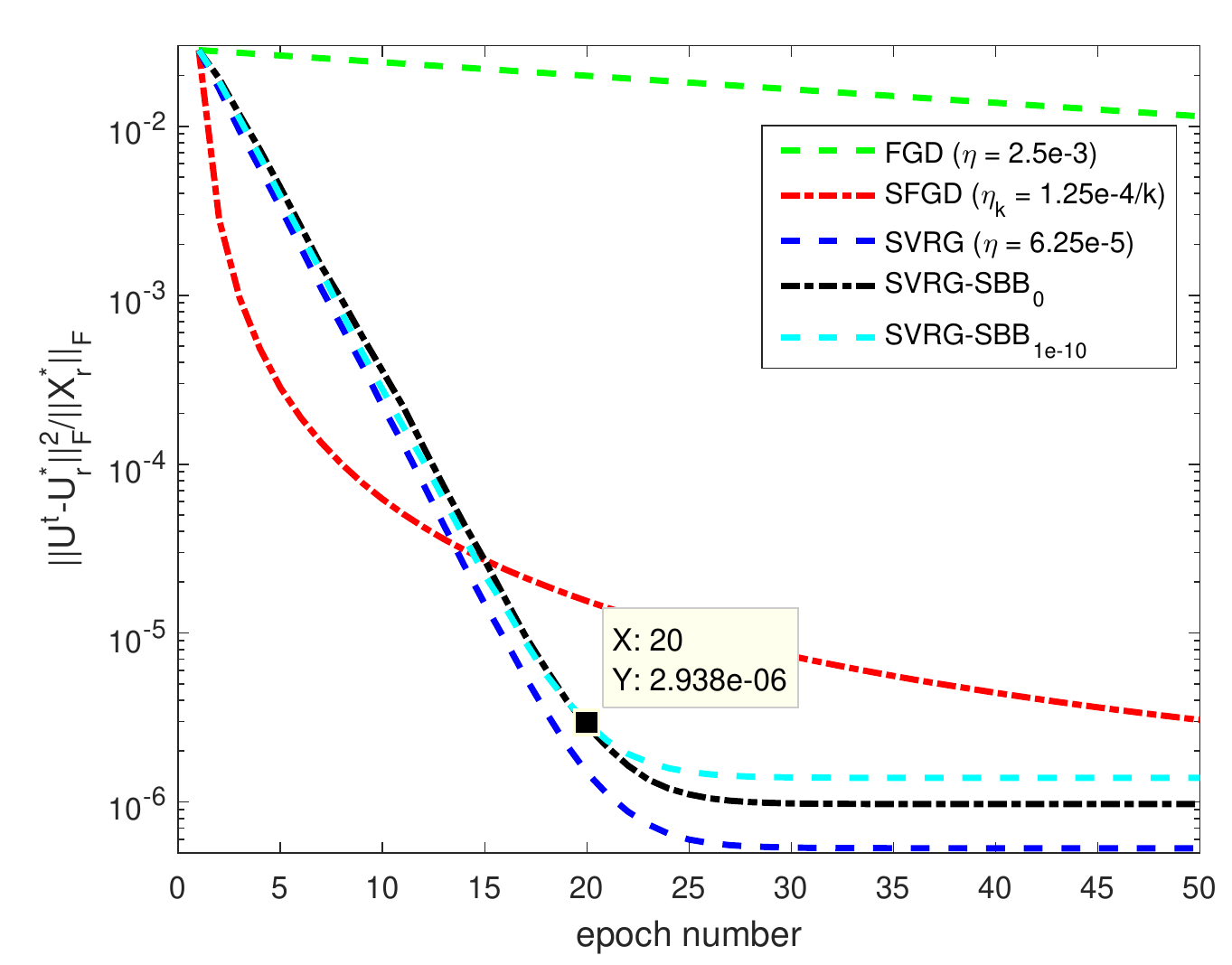}
\includegraphics[width=.47\linewidth]{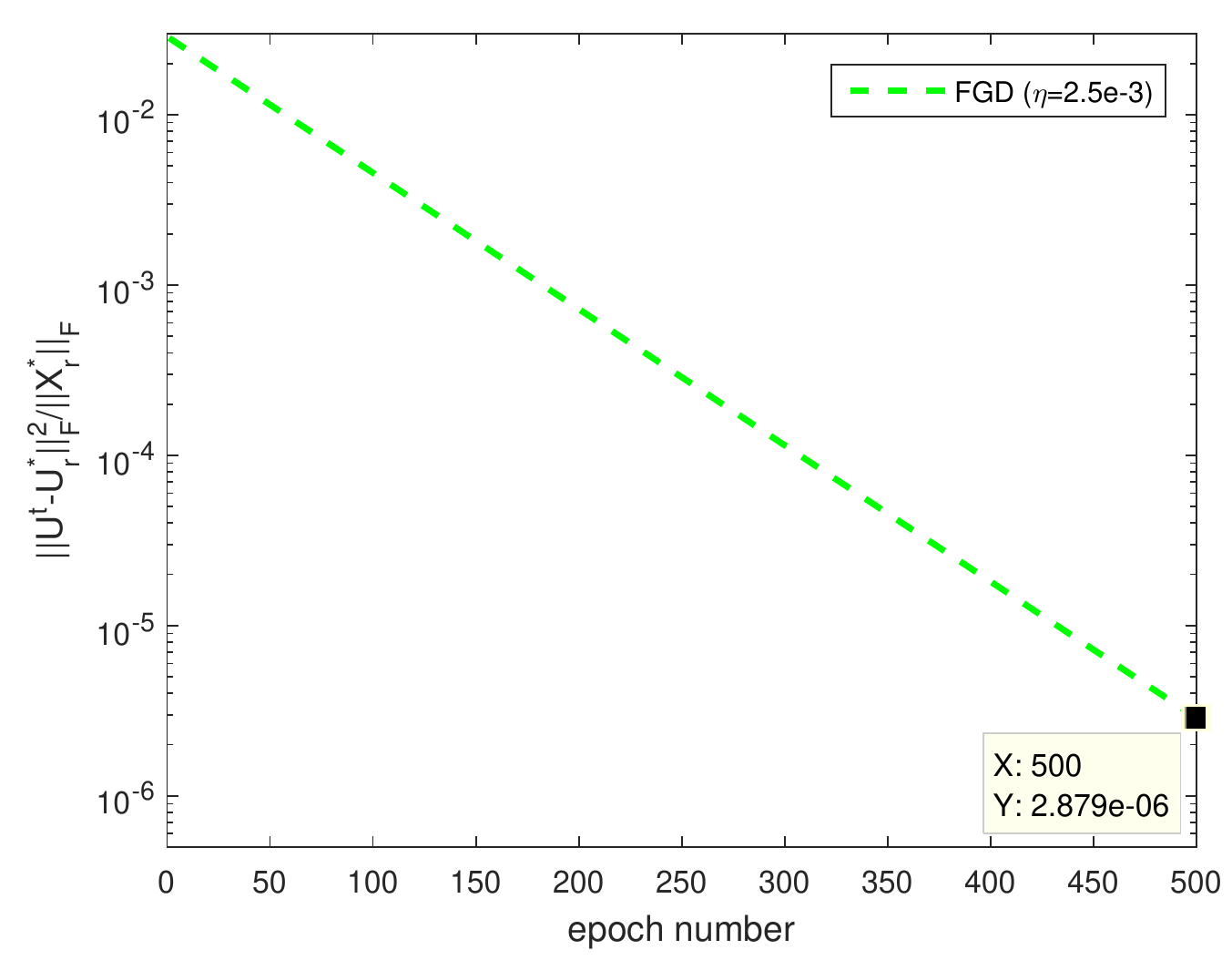}
%\vspace{-.3in}
\caption{Experiments for matrix sensing problem.
Left: trends of iterative errors of five algorithms. Right: trend of iterative error of FGD.
It requires about 20, 50 and 500 epochs for SVRG, SFGD and FGD respectively, to achieve the precision $3\times 10^{-6}$.
}
\label{Fig:matrix_sensing}
\end{figure}

Specifically, we let $p=5000$, the optimal matrix $X^* : = U^*{U^*}^T$ be a low-rank matrix with $\mathrm{rank}(X^*) = 5$ and the sample size $n=10p$. In such high-dimensional regime, the generic semidefinite optimization methods generally do not work. Therefore, we only compare the performance of the low-rank factorization based methods, i.e., FGD \cite{Sanghavi-FGD2016}, SFGD, and SVRG with three different step sizes studied in this paper.
%For FGD, we use $\eta = \frac{1}{4L}$, where $L$ is the Lipschitz constant of $\nabla f$. For SFGD, we use a more consecutive fixed step size $\hat{\eta} = \frac{1}{8L\|X^*\|_F}$ and a diminishing step size $\eta_t = \frac{\hat{\eta}}{t+1}$.
In this experiment, $r$ is set as $r^*$, and the initialization is constructed via the optimum $U^*$ with a random perturbation, and the step sizes for all algorithms are tuned in the hand-optimal way (shown in the figure).
%\commyy{better say what are they in appendix?}.
For three SVRG algorithms, the update frequency of the inner loop $m $ is set as the sample size $n$.
The experiment results are shown in Figure \ref{Fig:matrix_sensing}.
An epoch of SFGD includes $n$ iterations of SFGD, an epoch of FGD is exactly an iteration of FGD, and an epoch of SVRG is an iteration of outer loop.
The iterative error curves of SVRG, SFGD and FGD are shown along epochs since all of them exploit a full scan of gradients over sample per epoch and their computational complexities per epoch are thus comparable.
%, as analyzed latter.

From Figure \ref{Fig:matrix_sensing}, we can observe that all three SVRG algorithms converge exponentially fast to the global optimum with high precisions. To achieve the precision $3\times 10^{-6}$, it requires about 50 and 500 epochs for SFGD and FGD, respectively, while about 20 epochs are generally sufficient for three SVRG algorithms. In terms of epoch number, the considered SVRG methods are much faster than both FGD and SFGD.
These experiment results demonstrate the effectiveness of SVRG and also verify our developed theoretical results.

\subsection{Ordinal Embedding}

In this subsection, we apply SVRG to the ordinal embedding problem,
of which the Stochastic Triplet Embedding (STE) \cite{vandermaaten2012stochastic} is one of the typical models. The objective function is shown as follows:
\begin{align*}
%\label{Eq:ordinalembedding}
f(X) = \frac{1}{|{\cal C}|}\sum_{c\in {\cal C}} \ell_c(X) + \lambda \cdot \mathrm{tr}(X),
\end{align*}
where ${\cal C}$ is a set of ordinal constraints, $|{\cal C}|$ is its cardinality, and $\ell_c$ is the logistic loss.
%In this application, we use the stabilized BB step size with $\epsilon =0$ (denoted as SBB$_0$) instead of the original BB step size since the objective function $f$ in this application might be nonconvex.
To show the effectiveness of the considered SVRG methods, we compare the performance of SVRG (using fixed, SBB$_0$ and SBB$_\epsilon$ step sizes, where $\epsilon = 0.02$) with SFGD, FGD and the projected gradient descent (ProjGD) method, where the last two are batch methods.

\textbf{A. Music artist dataset:} We implement SVRG on the first real world dataset called \textit{Music artist dataset}, collected by \cite{ellis2002quest} via a web-based survey. In this dataset, there are $1032$ users and $412$ music artists. The number of triplets on the similarity between music artists is $213472$. A \textit{triplet} $(i, j, k)$ indicates an ordinal constraint like $d^2_{ij}(X)\leq d^2_{ik}(X)$, which means that ``\textit{music artist $i$ is more similar to artist $j$ than artist $k$}'', where $d^2_{ij}(X)$ is the Euclidean distance between artists $i$ and $j$, $i, j, k \in \{1,\ldots,p\}$, and $p$ is the number of total kinds of music artists.
Specifically, we use the data pre-processed by \cite{vandermaaten2012stochastic} via removing the inconsistent triplets from the original dataset. In this dataset, there are $9107$ triplets for $p=400$ artists. The genre labels for all artists are gathered using Wikipedia,
to distinguish nine music genres (rock, metal, pop, dance, hip hop, jazz, country, gospel, and reggae).
%Thus, the desired embedding dimension $r^* = 9$.
%\commyy{we set our desired embedding dimension $r=9$? ($r^*$ is unknown...)}.

\begin{figure}[h]
%\vspace{-.3in}
%\centerline{\fbox{This figure intentionally left non-blank}}
\includegraphics[width=.47\linewidth]{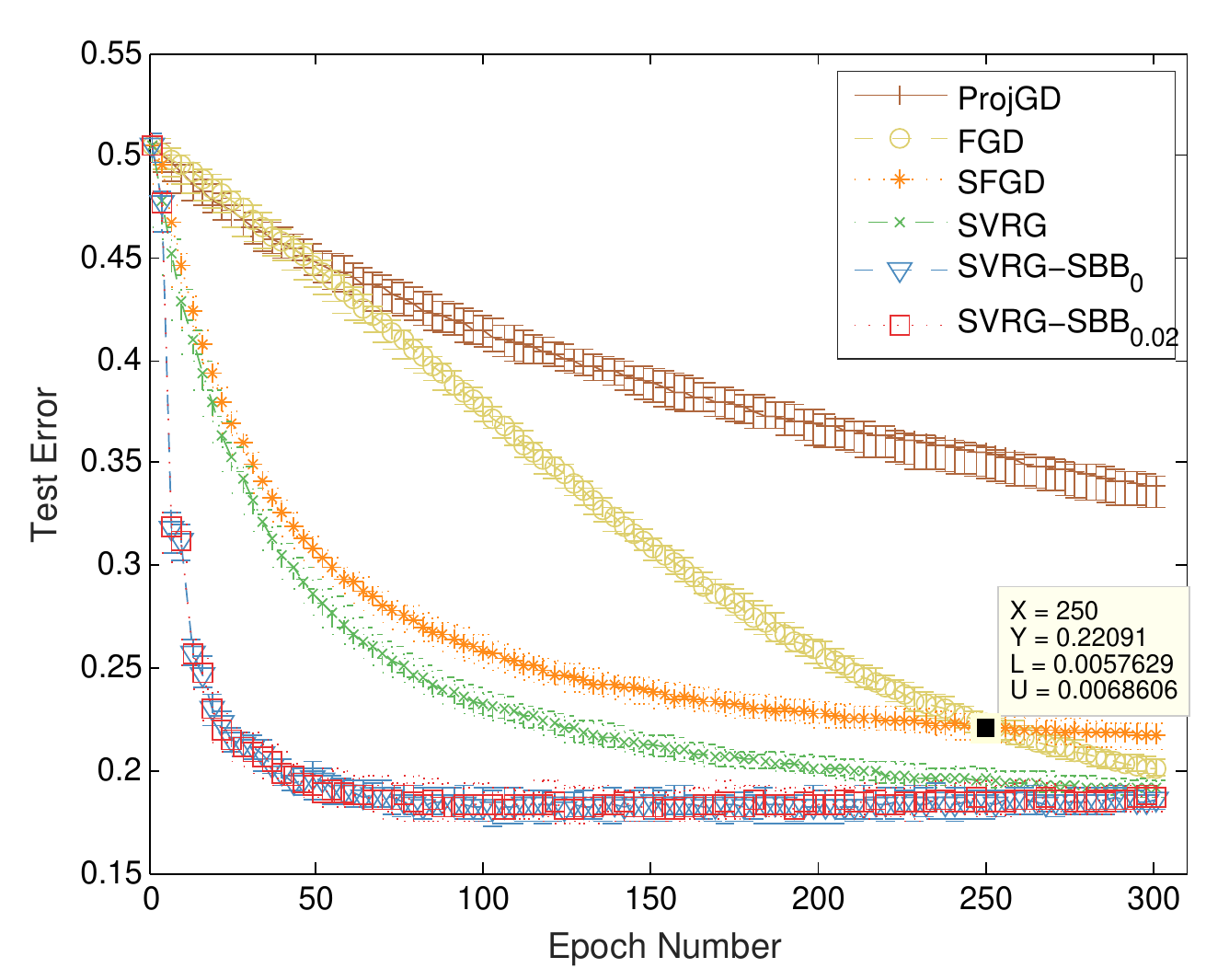}
\includegraphics[width=.47\linewidth]{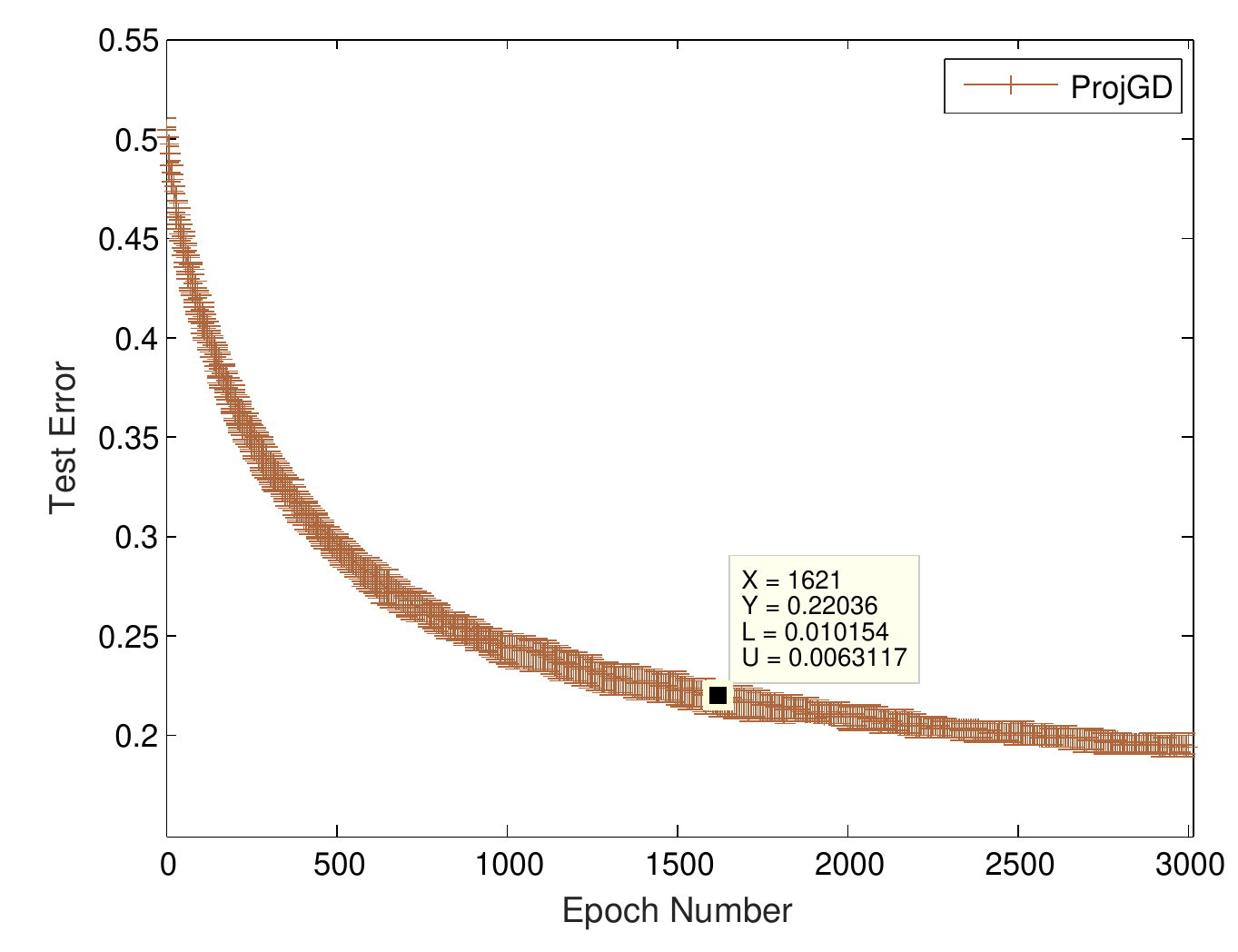}
%\vspace{-.3in}
\caption{Experiments for Music artist data.
To achieve the test error $0.22$,
about 40 epochs for SVRG-SBB$_0$ and SVRG-SBB$_{0.02}$, and 130 epochs for SVRG (fixed step size), and 260 epochs for both SFGD and FGD, and 1600 epochs for ProjGD are required.
}
\label{Fig:ordinal_embedding}
\end{figure}

For each method, we implement independently 50 trials, and then record their test errors. For each trail, $80\%$ triplets are randomly picked as the training set and the rest as the test set. All methods start with the same initial point, which is chosen randomly.
Each curve in Figure \ref{Fig:ordinal_embedding} shows the trend of test error of one method with respect to the  epoch number.

From Figure \ref{Fig:ordinal_embedding},
SVRG with SBB step sizes can significantly speed up SFGD and the batch methods in terms of epoch number. Particularly, the test error curves of two SVRG-SBB methods decay much faster than those of SFGD, FGD and ProjGD at the initial 50 epochs.

%In this subsection, we apply SVRG to the ordinal embedding problem, of which the generalized non-metric multidimensional scaling (GNMDS) \cite{Agarwal2007} is one of the typical models. The objective function is shown as follows:
%\begin{align*}
%%\label{Eq:ordinalembedding}
%f(X) = \frac{1}{|{\cal C}|}\sum_{c\in {\cal C}} \ell_c(X) + \lambda \cdot \mathrm{tr}(X),
%\end{align*}
%where ${\cal C}$ is a set of ordinal constraints, $|{\cal C}|$ is its cardinality, and $\ell_c$ is some loss function.
%In this application, we use the stabilized BB step size with $\epsilon =0$ (denoted as SBB$_0$) instead of the original BB step size since the objective function $f$ in this application might be nonconvex (see, \cite{Ma-SVRG-SBB2017} for detail). To show the effectiveness of the considered SVRG methods, we compare the performance of SVRG (with fixed, SBB$_0$ and SBB$_\epsilon$ step sizes, where $\epsilon = 0.1$) with SFGD, FGD and the projected gradient descent (ProjGD) method, where the last two are batch methods.

\textbf{B. eurodist dataset:}
We implement SVRG on another real world dataset called \textit{eurodist dataset}, which describes the ``driving'' distances between 21 cities in Europe, and is available in the stats library of {\bf R}. In this dataset, there are 21945 comparisons in total. A \textit{quadruplet} $(i, j, k, l)$ indicates an ordinal constraint like $d^2_{ij}(X)\leq d^2_{kl}(X)$, which means that ``\textit{the distance between cities $i$ and $j$ is shorter than the distance between cities $l$ and $k$}'', where $d^2_{ij}(X)$ is the ``driving'' distance between cities $i$ and $j$, $i, j, k, l \in \{1,\ldots,21\}$. One of the main task of this dataset is to embed these 21 cities in 2-dimensional space.
%Thus, the desired embedding dimension $r^* = 2$.

\begin{figure}[!t]
	\begin{minipage}[b]{0.95\linewidth}
		\centering
		%\vspace{-3.8cm}
		\includegraphics*[scale=0.4]{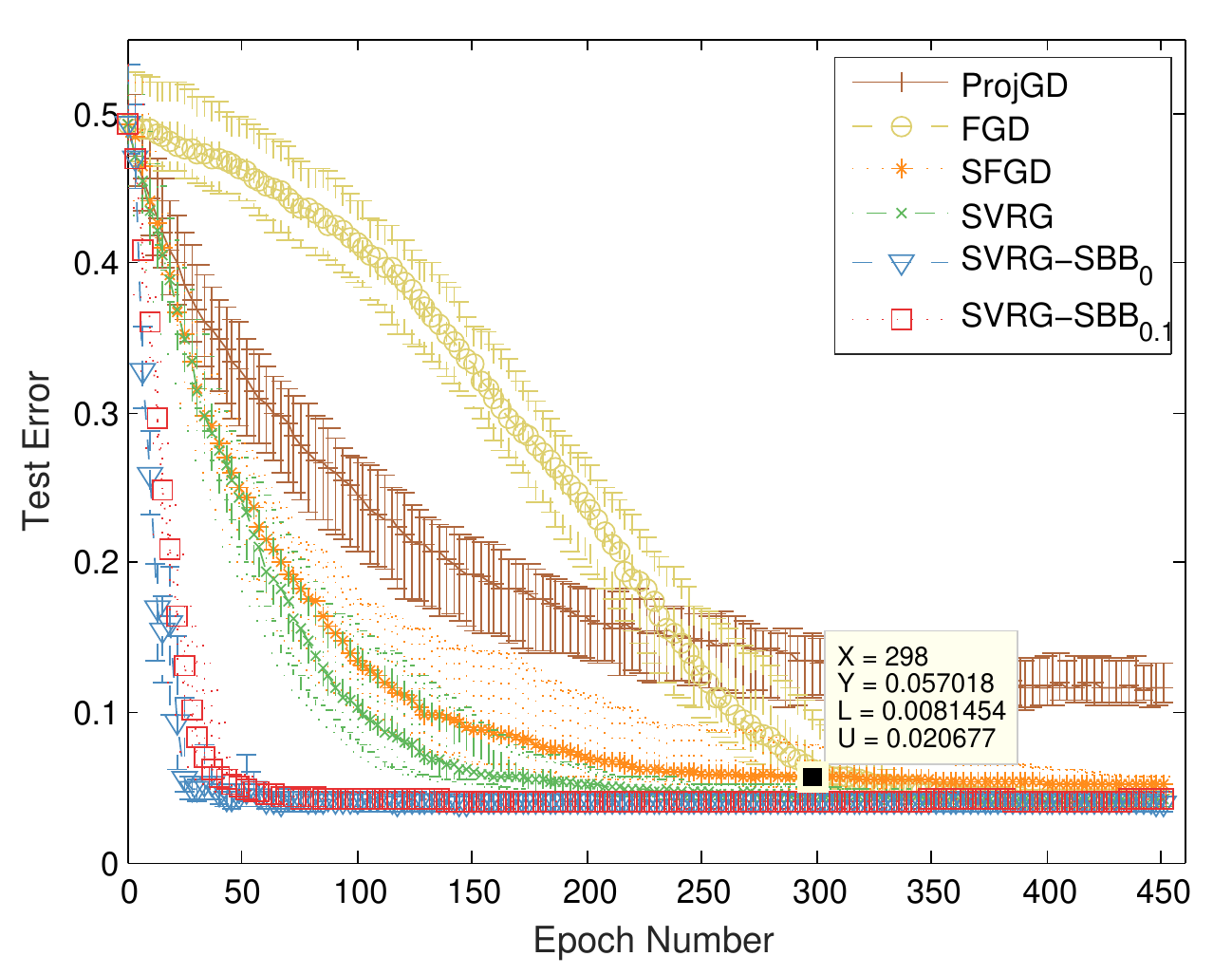}
		%\vspace{-4cm}
		%\centerline{{\small (a) Iterative error}}
	\end{minipage}
	\hfill
	\caption{Experiments for \textit{eurodist} dataset. To achieve the test error $0.057$,
about 40 epochs for SVRG-SBB$_0$ and SVRG-SBB$_{0.1}$, and 150 epochs for SVRG with a fixed step size, and 300 epochs for SFGD and FGD, and more than 450 epochs for ProjGD are required.
}
	\label{Fig:ordinal_embedding_eurodist}
\end{figure}

In this experiment, we first abstract all 3990 triplet ordinal comparisons from the total data set, and then use these triplets for learning. A \textit{triplet} $(i, j, k)$\footnote{A triplet $(i,j,k)$ is a special quadruplet $(i,k,j,k)$. We only use all triplets in this experiment because the existing and our codes are only suitable for dealing with triplet ordinal constraints. We will further prepare the codes for dealing with the quadruplet ordinal constraints.} indicates an ordinal constraint like $d^2_{ik}(X)\leq d^2_{jk}(X)$, which means that ``\textit{the distance between cities $i$ and $k$ is less than the distance between cities $j$ and $k$}''.
For each method, we implement independently 50 trials, and then record their test errors. For each trail, $80\%$ triplets are randomly picked as the training set and the rest as the test set. All methods start with the same initial point, which is chosen randomly.
Each curve in Figure \ref{Fig:ordinal_embedding_eurodist} shows the trend of test error of one method with respect to the  epoch number.
%The embedding results of all algorithms via comparing with the classical MDS (using the Matlab commond: mds.m) which
%is allowed to use the actual distance scores between all cities, are shown in (Supplementary Material: Appendix B).

From Figure \ref{Fig:ordinal_embedding_eurodist},
SVRG with SBB step sizes can speed up SFGD and both batch methods in terms of epoch number. Particularly, the test error curves of two SVRG-SBB methods decay much faster than those of SFGD, FGD and ProjGD at the initial 50 epochs.

%\begin{figure}[!t]
%	\begin{minipage}[b]{0.95\linewidth}
%		\centering
%		%\vspace{-3.8cm}
%		\includegraphics*[scale=0.4]{music_gnmds_testerror_crop.pdf}
%		%\vspace{-4cm}
%		%\centerline{{\small (a) Iterative error}}
%	\end{minipage}
%	\hfill
%	\caption{Experiments for ordinal embedding. To achieve the test error $0.2$,
%about 40 epochs for SVRG-SBB$_0$ and SVRG-SBB$_{0.025}$, and 150 epochs for SVRG (fixed step size), and 200 epochs for both SFGD and ProjGD, and 300 epochs for FGD are required.
%%Thus, there are at least three times of speedup of SVRG-SBB over the other methods.
%}
%	\label{Fig:ordinal_embedding}
%\end{figure}

\section{Conclusion}
\label{sc:conclusion}

In this paper, we consider a nonconvex stochastic semidefinite optimization problem, which emerges in many fields of science and engineering.
For the first time up to our knowledge, provable global linear convergence is established for stochastic variance reduced gradient (SVRG) algorithms to solve this nonconvex problem. Specifically, under common assumptions of restricted strong convexity of the objective function and small rank-$r$ approximation error, we can show that SVRG can converge to a global optimum at a linear rate as long as the initialization lies in a neighborhood of the optimum. The initial choice condition significantly improves the existing results for deterministic gradient descent. Moreover, our choice of step sizes includes both fixed and adaptive ones using Barzilai-Borwein (BB) step size with stabilization in nonconvex settings. Application examples show that the proposed scheme is promising in fast solving some large scale problems.

\section*{Acknowledgment}
The work of Jinshan Zeng is supported in part by the National Natural Science Foundation (NNSF) of China (No.61603162, 11501440), and the Doctoral start-up foundation of Jiangxi Normal University.
Yuan Yao's work is supported in part by HKRGC grant 16303817, 973 Program of China (No. 2015CB85600, 2012CB825501), NNSF of China (No. 61370004, 11421110001), as well as grants from Tencent AI Lab, Si Family Foundation, Baidu BDI, and Microsoft Research-Asia.
%The work of Ke Ma is supported by National Key Research and Development Plan (No.2016YFB0800403), National Natural Science Foundation of China (No.U1605252 and 61733007).

\newpage
\section*{Supplementary Material: Proofs}
\label{sc:proof}

For any matrix $U \in \mathbb{R}^{p\times r}$, let $Q_U$ be a basis of the column space of $U$. Denote ${\cal P}_U := Q_UQ_U^T$. Then ${\cal P}_U \cdot U = U$. For any matrix $Y \in \mathbb{R}^{p\times p}$, ${\cal P}_UY$ is a projection of $Y$ onto the subspace spanned by $X:= UU^T$.

\subsection*{A. Proof of Lemma \ref{Lemm:innerloop}}

%In the following, we describe a key lemma for the convergence of SVRG.

%\begin{Lemm}[A key lemma]
%\label{Lemm:innerloop}
%Let $\{U^t\}_{t=0}^m$ be the sequence at the $k$-th inner loop. Let Assumptions \ref{Assump:objfun}, \ref{Assump:r-approx error} and \ref{Assump:moment} hold. Let $\eta_k \in (0,\eta_{\max})$. If $\gamma_l<\mathbb{E}[\|\tilde{U}^k - U_r^*\|_F^2] < \gamma_u$, then the sequence
%$\{\mathbb{E}[\|U^t - U_r^*\|_F^2]\}$ is monotonically decreasing for $t=0,\ldots,m-1$, and
%\begin{align}
%\label{Eq:recur-innerloop}
%&\mathbb{E}_{i_t}[\|U^{t+1} - U_r^*\|_F^2] \\
%&\leq \|U^t - U_r^*\|_F^2 - \frac{2(\sqrt{2}-1)}{3}\eta_k \mu \sigma_r(X_r^*) \|U^t - U_r^*\|_F^2 \nonumber\\
%& + \frac{\eta_k L}{2\xi} \|U^t - U_r^*\|_F^4 + \eta_k^2 B_2 \cdot \|\tilde{U}^k - U_r^*\|_F^2 \nonumber\\
%&+ \frac{\eta_k L}{2} \|X^* - X_r^*\|_F^2,
%\end{align}
%where $B_2$ is specified in \eqref{Eq:B2};
%while if $\mathbb{E}[\|\tilde{U}^k - U_r^*\|_F^2]\leq \gamma_l$, then $\mathbb{E}[\|U^t - U_r^*\|_F^2]\leq \gamma_l$ for any $t=0,\ldots,m-1.$
%\end{Lemm}

The sketch proof of Lemma \ref{Lemm:innerloop} is show as follows.
We prove this lemma by induction. Specifically, we first show that if $\gamma_l < \mathbb{E}[\|U^t - U_r^*\|_F^2] < \gamma_u$, then $\mathbb{E}[\|U^{t+1} - U_r^*\|_F^2]\leq \mathbb{E}[\|U^t - U_r^*\|_F^2]  < \gamma_u$ for $t=0,\ldots,m-1$. Furthermore, $\mathbb{E}[\|U^{t+1} - U_r^*\|_F^2]$ can be estimated via noting that
\begin{align*}
&\mathbb{E}_{i_t}[\|U^{t+1} - U_r^*\|_F^2] \\
&= \|U^t - U_r^*\|_F^2 + \eta_k^2 \mathbb{E}_{i_t}[\|v_k^t\|_F^2] \\
&- 2\eta_k \mathbb{E}_{i_t}[\langle v_k^t, U^t - U_r^* \rangle],
\end{align*}
where $v_k^t = \nabla f_{i_t}(X^t)U^t - \nabla f_{i_t}(\tilde{X}^k)\tilde{U}^k + \nabla f(\tilde{X}^k)\tilde{U}^k$,
and establishing the bounds of both $\mathbb{E}_{i_t}[\|v_k^t\|_F^2]$ and $\mathbb{E}_{i_t}[\langle v_k^t, U^t - U_r^* \rangle]$ shown as the following two lemmas, respectively.

\begin{Lemm}[Bound of ${2\mathbb{E}_{i_t}[\langle v_k^t, U^t-U_r^* \rangle]}$]
\label{Lemm:inner-product}
Let Assumptions \ref{Assump:objfun} and \ref{Assump:moment} hold. Let $\{U^t\}_{t=0}^{m-1}$ be a sequence generated by SVRG in Algorithm \ref{alg:SVRG} at the $k$-th inner loop. Let $X^t = U^t{U^t}^T$ and $v_k^t = \nabla f_{i_t}(X^t)U^t - \nabla f_{i_t}(\tilde{X}^k)\tilde{U}^k + \nabla f(\tilde{X}^k)\tilde{U}^k$. If $\|U^t-U_r^*\|_F^2 <\gamma_0\sigma_r(X_r^*)$, then there holds
\begin{align*}
&2\mathbb{E}_{i_t}[\langle v_k^t, U^t-U_r^* \rangle]\\
&\geq \frac{\mu}{2}\|X^t - X_r^*\|_F^2 + \frac{\xi}{2L} \|{\cal P}_{U^t}\nabla f(X^t)\|_F^2\\
&- \frac{L}{2}\|X^* - X_r^*\|_F^2 - \frac{L}{2\xi}\|U^t-U_r^*\|_F^4.
\end{align*}
where $\xi$ is specified in \eqref{xi}.
\end{Lemm}

\textbf{Proof.}
%\begin{proof}
By Assumption \ref{Assump:moment},
\begin{align}
&2\mathbb{E}_{i_t}[\langle v_k^t, U^t-U_r^*\rangle]\nonumber\\
& =2 \langle \nabla f(X^t)U^t, U^t-U_r^* \rangle \nonumber\\
& = 2 \langle \nabla f(X^t), X^t-U_r^* {U^t}^T \rangle\nonumber\\
& = \langle \nabla f(X^t), X^t-X_r^* \rangle \nonumber\\
& + \langle \nabla f(X^t), X^t + X_r^* - 2 U_r^* {U^t}^T \rangle. \label{Eq:lowerbound}
%& = \langle \nabla f(X), X-X_r^* \rangle + \langle \nabla f(X),(U-V_U^*)(U-V_U^*)^T \rangle,
\end{align}
To bound the first term of \eqref{Eq:lowerbound}, we utilize the following three inequalities mainly by the Lipschitz differentiability and restricted strong convexity of $f$, that is,
%\begin{align*}
%\text{(i)}\quad &f(X_r^*) \geq f(X^t) + \langle \nabla f(X^t), X_r^*-X^t \rangle + \frac{\mu}{2}\|X_r^* - X^t\|_F^2, \\
%\text{(ii)}\quad & f(X^t) \geq f(X^*) + (1-\bar{\eta}/2)\bar{\eta}L^{-1}\cdot\|{\cal P}_{U^t}\nabla f(X^t)\|_F^2, \\
%\text{(iii)}\quad & f(X^*) \geq f(X_r^*) - \frac{L}{2}\|X^* - X_r^*\|_F^2,
%\end{align*}
\begin{align*}
\text{(i)} \  &f(X_r^*) \geq \\
\quad & f(X^t) + \langle \nabla f(X^t), X_r^*-X^t \rangle + \frac{\mu}{2}\|X_r^* - X^t\|_F^2, \\
\text{(ii)} \   & f(X^t) \geq f(X^*) + (1-\bar{\eta}/2)\bar{\eta}L^{-1}\cdot\|{\cal P}_{U^t}\nabla f(X^t)\|_F^2, \\
\text{(iii)}\   & f(X^*) \geq f(X_r^*) - \frac{L}{2}\|X^* - X_r^*\|_F^2,
\end{align*}
where (i) holds for the $(\mu,r)$-restricted strong convexity of $f$, (ii) holds the following inequality induced by the $L$-Lipschitz differentiability of $f$, i.e.,
\begin{align*}
f(X^t)
&\geq f(\bar{X}) + \langle \nabla f(X^t), X^t-\bar{X} \rangle - \frac{L}{2}\|X^t-\bar{X}\|_F^2 \\
&(\text{where}\ \bar{X}: = X^t - \frac{\bar{\eta}}{L}{\cal P}_{U^t}\nabla f(X^t){\cal P}_{U^t} ) \\
&= f(\bar{X}) +(1-\bar{\eta}/2)\bar{\eta}L^{-1} \cdot \|{\cal P}_{U^t}\nabla f(X^t)\|_F^2,
\end{align*}
and $f(\bar{X}) \geq f(X^*)$ since $X^*$ is an optimum and $\bar{X}$ is a feasible point by Lemma \ref{Lemm:feasibility}(b),
and (iii) holds for the $L$-Lipschitz differentiability of $f$ and the optimality condition $\nabla f(X^*)U^* =0$, i.e.,
\begin{align*}
f(X_r^*)
&\leq f(X^*) + \langle \nabla f(X^*), X^* - X_r^* \rangle + \frac{L}{2} \|X^* - X_r^*\|_F^2 \\
&= f(X^*) + \frac{L}{2}\|X^* - X_r^*\|_F^2,
\end{align*}
where the equality holds for $\nabla f(X^*)U^* =0$, which directly implies the following facts: $\nabla f(X^*)U_r^* =0$, $\nabla f(X^*)X^* = 0$ and $\nabla f(X^*)X_r^* =0$ due to $X^* = U^* {U^*}^T$ and $X_r^* = U_r^* {U_r^*}^T$. Summing the inequalities (i)-(iii) yields
\begin{align}
&\langle \nabla f(X^t), X^t-X_r^* \rangle \geq \frac{\mu}{2} \|X^t-X_r^*\|_F^2 \label{Eq:T11}\\
&+ (1-\bar{\eta}/2)\bar{\eta}L^{-1}\cdot \|{\cal P}_{U^t}\nabla f(X^t)\|_F^2 - \frac{L}{2}\|X^* - X_r^*\|_F^2. \nonumber
%\\
%&\geq (\sqrt{2}-1)\mu \sigma_r(X_r^*) {\cal E}(U,U_r^*) + (1-\bar{\eta}/2)\bar{\eta}L^{-1} \cdot \|{\cal P}_U \nabla f(X)\|_F^2 - \frac{L}{2}\|X^* - X_r^*\|_F^2, \nonumber
\end{align}
%where the second inequality is due to Lemma \ref{Lemm:dist-twomat}, i.e., $\|X^t-X_r^*\|_F^2 \geq 2(\sqrt{2}-1)\sigma_r(X_r^*) {\cal E}(U,U_r^*)$.
On the other hand, we observe that
\begin{align}
& \langle \nabla f(X^t), X^t + X_r^* - 2 U_r^* {U^t}^T \rangle \nonumber\\
& = \langle {\cal P}_{U^t} \nabla f(X^t) + ({\bf I}-{\cal P}_{U^t})\nabla f(X^t),  X^t + X_r^* - 2 U_r^* {U^t}^T \rangle \nonumber\\
& = \langle {\cal P}_{U^t} \nabla f(X^t), X^t + X_r^* - 2 U_r^* {U^t}^T \rangle \nonumber\\
& =\langle {\cal P}_{U^t} \nabla f(X^t), (U^t-U_r^*)(U-U_r^*)^T \rangle \nonumber\\
& \geq -\frac{(1-\bar{\eta}/2)\bar{\eta}}{2L} \|{\cal P}_{U^t}\nabla f(X^t)\|_F^2 \nonumber\\
& - \frac{L}{2\bar{\eta}(1-\bar{\eta}/2)} \cdot \|U^t - U_r^*\|_F^4. \label{Eq:T12}
\end{align}
where the second equality is due to $\langle ({\bf I}-{\cal P}_{U^t})\nabla f(X^t), X^t \rangle =0$, $\langle ({\bf I}-{\cal P}_{U^t})\nabla f(X^t), U_r^* {U^t}^T \rangle =0$ and $\langle ({\bf I}-{\cal P}_{U^t})\nabla f(X^t), X_r^* \rangle =0$ by Lemma \ref{Lemm:feasibility}(c),
the last equality holds for $X_r^* = U_r^* {U_r^*}^T$, and the inequality holds for the basic inequality: $\langle Y, Z \rangle \geq -\frac{c}{2} \|Y\|_F^2 - \frac{1}{2 c} \|Z\|_F^2$ for any $Y, Z \in \mathbb{R}^{p\times p}$ and $c>0$.
Substituting \eqref{Eq:T11} and \eqref{Eq:T12} into \eqref{Eq:lowerbound} concludes this lemma.
%\end{proof}
$\Box$

\begin{Lemm}[Bound of ${\mathbb{E}_{i_t}[\|v_k^t\|_F^2]}$]
\label{Lemm:vkt}
%Consider the $k$-th inner loop. Let $v_k^t : = \nabla f_{i_t}(X^t)U^t - \nabla f_{i_t}(\tilde{X}^k)\tilde{U}^k + g_k$.
Let Assumptions \ref{Assump:objfun}, \ref{Assump:r-approx error} and \ref{Assump:moment} hold. Assume that $\|U^t - U_r^*\|^2_F < \gamma_u$ and $\|\tilde{U}^k - U_r^*\|^2_F < \gamma_u$,
then
\begin{align*}
\mathbb{E}_{i_t}[\|v_k^t\|_F^2]
& \leq 4(B_0 + B_1)(\|U^t - U_r^*\|_F^2+\|\tilde{U}^k - U_r^*\|_F^2) \\
& + 4 L^2 {\cal B}(\|X^t - X_r^*\|_F^2 + \|\tilde{X}^k - X_r^*\|_F^2) \\
& + \|{\cal P}_{U^t} \nabla f(X^t)\|_F^2 \cdot \|X^t\|_F.
\end{align*}
\end{Lemm}

\textbf{Proof.}
%\begin{proof}
Note that
\begin{align*}
\|v_k^t\|_F^2
& = \|\nabla f_{i_t}(X^t)U^t - \nabla f_{i_t}(\tilde{X}^k) \tilde{U}^k\|_F^2 \\
& + \|\nabla f(\tilde{X}^k)\tilde{U}^k\|_F^2 \\
& + 2\langle \nabla f_{i_t}(X^t)U^t - \nabla f_{i_t}(\tilde{X}^k)\tilde{U}^k, \nabla f(\tilde{X}^k)\tilde{U}^k \rangle.
\end{align*}
Thus,
\begin{align}
\label{Eq:E-vkt}
&\mathbb{E}_{i_t}[\|v_k^t\|_F^2]\nonumber\\
&= \mathbb{E}_{i_t}[\|\nabla f_{i_t}(X^t)U^t - \nabla f_{i_t}(\tilde{X}^k) \tilde{U}^k\|_F^2] \nonumber\\
&+\|\nabla f(\tilde{X}^k)\tilde{U}^k\|_F^2\nonumber\\
&+ 2\langle \nabla f(X^t)U^t - \nabla f(\tilde{X}^k)\tilde{U}^k, \nabla f(\tilde{X}^k)\tilde{U}^k \rangle\nonumber\\
& = \mathbb{E}_{i_t}[\|\nabla f_{i_t}(X^t)U^t - \nabla f_{i_t}(\tilde{X}^k) \tilde{U}^k\|_F^2] \nonumber\\
&- \|\nabla f(X^t)U^t - \nabla f(\tilde{X}^k)\tilde{U}^k\|_F^2 + \|\nabla f(X^t)U^t\|_F^2 \nonumber\\
& \leq \mathbb{E}_{i_t}[\|\nabla f_{i_t}(X^t)U^t - \nabla f_{i_t}(\tilde{X}^k) \tilde{U}^k\|_F^2] \nonumber\\
&+ \|{\cal P}_{U^t} \nabla f(X^t)\|_F^2 \cdot \|X^t\|_F,
\end{align}
where the last inequality holds for
$\|\nabla f(X^t)U^t - \nabla f(\tilde{X}^k)\tilde{U}^k\|_F^2 \geq 0$
and
\begin{align*}
&\|\nabla f(X^t)U^t\|_F^2 \\
&= \|{\cal P}_{U^t} \nabla f(X^t) U^t + ({\bf I}-{\cal P}_{U^t}) \nabla f(X^t)U^t\|_F^2 \\
&= \|{\cal P}_{U^t} \nabla f(X^t) U^t\|^2 \\
&\leq \|{\cal P}_{U^t} \nabla f(X^t)\|_F^2 \cdot \|X^t\|_F.\\
\end{align*}
In the following, we bound the first term of \eqref{Eq:E-vkt}. Note that
\begin{align*}
&\|\nabla f_{i_t}(X^t)U^t - \nabla f_{i_t}(\tilde{X}^k) \tilde{U}^k\|_F^2\\
&=\|\nabla f_{i_t}(X^t)(U^t-\tilde{U}^k) + (\nabla f_{i_t}(X^t)-\nabla f_{i_t}(\tilde{X}^k)) \tilde{U}^k\|_F^2 \\
&\leq 2\|\nabla f_{i_t}(X^t)(U^t-\tilde{U}^k)\|_F^2 \\
& + 2 \|(\nabla f_{i_t}(X^t)-\nabla f_{i_t}(\tilde{X}^k)) \tilde{U}^k\|_F^2\\
&\leq 2\|\nabla f_{i_t}(X^t)\|_F^2 \|U^t-\tilde{U}^k\|_F^2 \\
& + 2L^2 \|\tilde{X}^k\|_F \|X^t - \tilde{X}^k\|_F^2,
\end{align*}
which follows
\begin{align*}
&\mathbb{E}_{i_t}[\|\nabla f_{i_t}(X^t)U^t - \nabla f_{i_t}(\tilde{X}^k) \tilde{U}^k\|_F^2]\\
& \leq  2(\mathbb{E}_{i_t}[\|\nabla f_{i_t}(X^t)\|_F^2] - \|\nabla f(X^t)\|_F^2)\|U^t-\tilde{U}^k\|_F^2 \\
&+ 2 \|\nabla f(X^t)\|_F^2 \cdot \|U^t-\tilde{U}^k\|_F^2 + 2L^2 \|\tilde{X}^k\|_F \|X^t - \tilde{X}^k\|_F^2\\
&\leq 2(B_0 + B_1)\cdot \|U^t-\tilde{U}^k\|_F^2 + 2L^2 {\cal B}\|X^t - \tilde{X}^k\|_F^2\\
&\leq 4(B_0 + B_1)\cdot (\|U^t-U_r^*\|_F^2 + \|\tilde{U}^k-U_r^*\|_F^2) \\
&+ 4L^2 {\cal B} (\|X^t - X_r^*\|_F^2 + \|\tilde{X}^k - X_r^*\|_F^2),
\end{align*}
where $B_0$, $B_1$ and ${\cal B}$ are specified before \eqref{Eq:B2}.
%in \eqref{Eq:B0}, \eqref{Eq:B1} and \eqref{Eq:B-cal}, respectively.
%\begin{align*}
%&B_0 := \sup_{\{U:\|U-U_r^*\|_F < \sqrt{\gamma_0} \sigma_r(U_r^*)\}} (\mathbb{E}_{i_t}[\|\nabla f_{i_t}(UU^T)\|_F^2] - \|\nabla f(UU^T)\|_F^2),\\
%&B_1 := \sup_{\{U:\|U-U_r^*\|_F < \sqrt{\gamma_0} \sigma_r(U_r^*)\}} \|\nabla f(UU^T)\|_F,\\
%&{\cal B}:= \sup_{\{U:\|U-U_r^*\|_F < \sqrt{\gamma_0} \sigma_r(U_r^*)\}} \|UU^T\|_F.
%\end{align*}
Substituting the above inequality into \eqref{Eq:E-vkt}, we can conclude this lemma.
%\end{proof}
$\Box$

Based on the above two lemmas, we give the proof of Lemma \ref{Lemm:innerloop}.

\textbf{Proof of Lemma \ref{Lemm:innerloop}:}
%Taking expectation of both sides of the above equality over $i_t$, it becomes
By Lemma \ref{Lemm:inner-product} and Lemma \ref{Lemm:vkt},
\begin{align*}
&\mathbb{E}_{i_t}[\|U^{t+1} - U_r^*\|_F^2] \leq \|U^t - U_r^*\|_F^2\nonumber\\
& - \eta_k \left[ \frac{\mu}{2}\|X^t - X_r^*\|_F^2 + \frac{\xi}{2L} \|{\cal P}_{U^t}\nabla f(X^t)\|_F^2 \right] \nonumber\\
& + \eta_k \left[ \frac{L}{2}\|X^* - X_r^*\|_F^2 + \frac{L}{2\xi}\|U^t-U_r^*\|_F^4\right] \nonumber\\
& + 4\eta_k^2 (B_0 + B_1)(\|U^t - U_r^*\|_F^2+\|\tilde{U}^k - U_r^*\|_F^2) \nonumber\\
& + 4\eta_k^2  L^2 {\cal B}\left(\|X^t - X_r^*\|_F^2 + \|\tilde{X}^k - X_r^*\|_F^2 \right)\nonumber\\
& + \eta_k^2 \|{\cal P}_{U^t} \nabla f(X^t)\|_F^2 \cdot \|X^t\|_F \nonumber\\
& \leq \|U^t - U_r^*\|_F^2 - \left( \frac{\eta_k \mu}{2} - 4\eta_k^2 L^2 {\cal B}\right)\|X^t - X_r^*\|_F^2 \nonumber\\
& + 4\eta_k^2 (B_0 + B_1)\|U^t - U_r^*\|_F^2 + \frac{\eta_k L}{2\xi} \|U^t - U_r^*\|_F^4 \nonumber\\
& + 4\eta_k^2 \left[ L^2 {\cal B} \|\tilde{X}^k - X_r^*\|_F^2 + (B_0 + B_1) \|\tilde{U}^k - U_r^*\|_F^2\right] \nonumber\\
& + \frac{\eta_k L}{2} \|X^* - X_r^*\|_F^2 \nonumber
\end{align*}
\begin{align*}
& \leq \|U^t - U_r^*\|_F^2 - \eta_k \|U^t - U_r^*\|_F^2 \times \nonumber\\
& \left[2(\sqrt{2}-1)\sigma_r(X_r^*)\left(\frac{\mu}{2}-4\eta_k L^2 {\cal B}\right) - 4\eta_k (B_0 + B_1) \right] \nonumber\\
& + \frac{\eta_k L}{2\xi} \|U^t - U_r^*\|_F^4 \nonumber\\
%2(\sqrt{2}-1)\eta\sigma_r(X_r^*)(\frac{\mu}{2}-4\eta \bar{L}^2 {\cal B})\|U^t - U_r^*\|_F^2 + 4\eta^2 (B_0 + B_1)\|U^t - U_r^*\|_F^2
& + 4\eta_k^2 \left[ L^2 {\cal B} \|\tilde{X}^k - X_r^*\|_F^2 + (B_0 + B_1) \|\tilde{U}^k - U_r^*\|_F^2\right] \nonumber\\
& + \frac{\eta_k L}{2} \|X^* - X_r^*\|_F^2,\nonumber\\
&\leq \|U^t - U_r^*\|_F^2 + \frac{\eta_k L}{2\xi} \|U^t - U_r^*\|_F^4- \eta_k \|U^t - U_r^*\|_F^2 \times \nonumber\\
& \left[2(\sqrt{2}-1)\sigma_r(X_r^*)\left(\frac{\mu}{2}-4\eta_k L^2 {\cal B}\right) - 4\eta_k (B_0 + B_1) \right]
\nonumber\\
& + 4\eta_k^2 \left[ 2L^2 {\cal B}({\cal B} + \|X_r^*\|_F) + B_0 + B_1 \right] \cdot \|\tilde{U}^k - U_r^*\|_F^2 \nonumber\\
& + \frac{\eta_k L}{2} \|X^* - X_r^*\|_F^2,
\end{align*}
where the first inequality is due to Lemma \ref{Lemm:inner-product} and Lemma \ref{Lemm:vkt},
the second inequality holds for  $\eta_k < \eta_{\max} \leq \frac{\xi}{2{\cal B}L}$,
the third inequality holds for $\eta_k < \eta_{\max} \leq \frac{1}{8{\cal B}\kappa L}$ and Lemma \ref{Lemm:dist-twomat}(b), the final inequality holds for Lemma \ref{Lemm:dist-twomat}(a).
Since $\eta_k < \eta_{\max} \leq \frac{(\sqrt{2}-1)\mu \sigma_r(X_r^*)}{12[2(\sqrt{2}-1)\sigma_r(X_r^*)L^2{\cal B} + B_0 + B_1]}$, then
\begin{align*}
& 2(\sqrt{2}-1)\sigma_r(X_r^*)\left(\frac{\mu}{2}-4\eta_k L^2 {\cal B}\right) - 4\eta_k (B_0 + B_1) \nonumber\\
& \geq \frac{2(\sqrt{2}-1)\mu \sigma_r(X_r^*)}{3}.
\end{align*}
Thus, substituting the above inequality into the first main inequality yields \eqref{Eq:recur-innerloop}.

Furthermore, by the assumption of this lemma and $\gamma_u \leq \gamma_0 \sigma_r(X_r^*)$, we have
\[\|\tilde{U}^k-U_r^*\|_F^2 < \gamma_0 \sigma_r(X_r^*).\]
Thus,
\begin{align}
&\mathbb{E}_{i_t}[\|U^{t+1} - U_r^*\|_F^2] \nonumber\\
&\leq \|U^t - U_r^*\|_F^2 + \frac{\eta_k L}{2\xi} \times\nonumber\\
& \left[\|U^t - U_r^*\|_F^4  + \xi \|X^* - X_r^*\|_F^2 + \frac{2 \xi \eta_k B_2}{L}\cdot \gamma_0 \sigma_r(X_r^*)\right] \nonumber\\
& - \frac{\eta_k L}{2\xi} \cdot \frac{4(\sqrt{2}-1)\xi \sigma_r(X_r^*)}{3\kappa} \|U^t - U_r^*\|_F^2
\label{Eq:Ut+1-Ur*}\\
%& + \frac{\eta_k L}{2\xi} \cdot \frac{8\eta_k \xi}{L} \left[ 2L^2 {\cal B}({\cal B}+\|X_r^*\|_F) + B_0 +B_1\right]\cdot \gamma_0 \sigma_r(X_r^*) \nonumber\\
& \leq \|U^t - U_r^*\|_F^2  +\frac{\eta_k L}{2\xi} \times \nonumber\\
& \left[\|U^t - U_r^*\|_F^4 - \frac{4(\sqrt{2}-1)\xi \sigma_r(X_r^*)}{3\kappa} \|U^t - U_r^*\|_F^2 \right] \nonumber\\
&+\frac{\eta_k L}{2\xi} \cdot \left[\xi \|X^* - X_r^*\|_F^2 + \frac{(\sqrt{2}-1)^2\xi^2 \sigma_r^2(X_r^*)}{9\kappa^2}\right], \nonumber
\end{align}
where the second inequality holds for $\eta_k < \eta_{\max} \leq \frac{(\sqrt{2}-1)\mu \xi \sigma_r(X_r^*)}{12 B_2}$, and $\eta_k < \eta_{\max} \leq \frac{(\sqrt{2}-1)^2 \xi \mu \sigma_r(X_r^*)}{18\kappa B_2 \gamma_0}$. By the definitions of $\gamma_l$ \eqref{gamma l} and $\gamma_u$ \eqref{gamma u}, the above inequality implies
\begin{align*}
&\mathbb{E}_{i_t}[\|U^{t+1} - U_r^*\|_F^2]
\leq \|U^t - U_r^*\|_F^2 \\
&- \frac{\eta_k L}{2\xi} (\gamma_u - \|U^t - U_r^*\|_F^2)(\|U^t - U_r^*\|_F^2 - \gamma_l),
\end{align*}
which implies $\mathbb{E}[\|U^{t+1} - U_r^*\|_F^2] \leq \mathbb{E}[\|U^{t} - U_r^*\|_F^2]$ if $\gamma_l<\mathbb{E}[\|U^{t} - U_r^*\|_F^2]<\gamma_u$.
Inductively, we can claim the first part of this lemma.

Define a univariate function $h(z) = z - \frac{\eta_k L}{2\xi} \cdot \left( (\gamma_l + \gamma_u) z - z^2 - \gamma_l \cdot \gamma_u \right)$ for any $z \in \mathbb{R}_+$. Then its derivative is
\begin{align*}
& h'(z) = 1- \frac{\eta_k L}{2\xi} \cdot (\gamma_l + \gamma_u) + \frac{\eta_k L}{\xi} \cdot z \\
& = 1-(\sqrt{2}-1)\eta_k \mu \sigma_r(X_r^*) + \frac{\eta_k L}{\xi} \cdot z >0,
\end{align*}
for $0<z\leq \gamma_l$,
where the second equality holds for \eqref{Eq:gamma-l-u}, and the inequality is due to $1-(\sqrt{2}-1)\eta_k \mu \sigma_r(X_r^*) >0$ for any $\eta_k \in (0,\eta_{\max})$.
 Thus, for any $0<z\leq \gamma_l$,
\[
h(z)\leq h(\gamma_l) = \gamma_l,
\]
which shows that the last statement of this lemma holds. Therefore, we end the proof of this lemma.
%\end{proof}
$\Box$

\subsection*{B. Proof of Corollary \ref{Coro:exact-recovery}}

\textbf{Proof.}
%\begin{proof}
Note that $\bar{\eta}_{\max} \leq \eta_{\max}$.
By Theorem \ref{Thm:svrg}, if
\begin{align*}
&\gamma_l:= \frac{(2-\sqrt{3})\cdot(\sqrt{2}-1)\xi\sigma_r(X_r^*)}{3\kappa}  \\
&<\|\tilde{U}^0-U_r^*\|_F^2
<\frac{(2+\sqrt{3})\cdot(\sqrt{2}-1)\xi\sigma_r(X_r^*)}{3\kappa} : = \gamma_u,
\end{align*}
then it is obvious that SVRG converges to the optimum $X^*$ at a linear rate. As a consequence, we only need to prove the exact recovery of SVRG when $\|\tilde{U}^0-U_r^*\|_F^2 \leq \gamma_l$.
By Theorem \ref{Thm:svrg}, in this case, $\mathbb{E}[\|\tilde{U}^k-U_r^*\|_F^2] \leq \gamma_l$ for all $k \in \mathbb{N}$. Actually, by the proof of Theorem \ref{Thm:svrg}, at any $k$-th inner loop,
\begin{align}
\label{Eq:smallbound}
\mathbb{E}[\|\tilde{U}^t-U_r^*\|_F^2] \leq \gamma_l
\end{align}
for any $t = 1,\ldots, m$.

In this case, it is obvious that Lemmas \ref{Lemm:inner-product} and \ref{Lemm:vkt} still hold, and \eqref{Eq:Ut+1-Ur*} in the proof of Lemma \ref{Lemm:innerloop} should be revised as
\begin{align}
&\mathbb{E}_{i_t}[\|U^{t+1} - U_r^*\|_F^2] \nonumber\\
&\leq \|U^t - U^*\|_F^2 +\frac{\eta_k L}{2\xi} \times \nonumber\\
& \left[\|U^t - U^*\|_F^4 - \frac{4(\sqrt{2}-1)\xi \sigma_r(X_r^*)}{3\kappa} \|U^t - U^*\|_F^2 \right]  \nonumber\\
&+ \eta_k^2 B_2 \cdot \|\tilde{U}^k - U^*\|_F^2 \nonumber\\
&\leq \left(1-\frac{\eta_k  L}{2\xi}\cdot \gamma_u \right) \cdot \|U^t - U^*\|_F^2 \nonumber\\
&+ \eta_k^2 B_2 \cdot \|\tilde{U}^k - U^*\|_F^2, \label{Eq:recur-innerloop-small}
\end{align}
where the second inequality holds for \eqref{Eq:smallbound} and \eqref{Eq:gamma-l-u}. By \eqref{Eq:recur-innerloop-small}, recursively, after some simplifications we have
\begin{align}
&\mathbb{E}[\|\tilde{U}^{k+1} - U^*\|_F^2] \label{Eq:Xk+1-Xk}\\
&\leq \left(1-\frac{\eta_k L}{2\xi}\cdot \gamma_u\right)^m \|\tilde{U}^k - U^*\|_F^2 \nonumber\\
& + \frac{2B_2\eta_k \xi}{L \gamma_u}  \cdot \left[ 1 - \left(1-\frac{\eta_k L}{2\xi}\cdot \gamma_u\right)^m\right] \|\tilde{U}^k - U^*\|_F^2 \nonumber
%& \leq \left(1-\frac{\eta_k L}{2\xi}\cdot \gamma_u\right)^m \mathbb{E}[\|\tilde{U}^k - U^*\|_F^2]
%+ \frac{16\eta_k L {\cal B}({\cal B} + \|X^*\|_F) \xi}{\gamma_u} \left[ 1 - \left(1-\frac{\eta_k L}{2\xi}\cdot \gamma_u\right)^m\right] \|\tilde{U}^k - U^*\|_F^2 \nonumber\\
%& = \left( \left(1-\frac{\eta_k L}{2\xi}\cdot \gamma_u\right)^m + \frac{2B_2\eta_k \xi}{L \gamma_u} \cdot \left[ 1 - \left(1-\frac{\eta_k L}{2\xi}\cdot \gamma_u\right)^m\right]\right)\times  \nonumber\\
%& \|\tilde{U}^k - U^*\|_F^2.
\end{align}
Since $\eta_k \in (0,\bar{\eta}_{\max})$, then
\[
0<\frac{2B_2 \eta_k \xi}{L \gamma_u} <1,
\]
and thus,
\begin{align*}
&0< \left(1-\frac{\eta_k L}{2\xi}\cdot \gamma_u\right)^m\\
& + \frac{2B_2 \eta_k \xi}{L \gamma_u}  \cdot \left[ 1 - \left(1-\frac{\eta_k L}{2\xi}\cdot \gamma_u\right)^m\right]
<1,
\end{align*}
which implies that SVRG converges to $X^*$ at a linear rate.
Therefore, we finish the proof of this corollary.
%\end{proof}
$\Box$

%\section*{Appendix}
%
%In the appendix, we first present several lemmas, which are frequently used in this paper, and then provide the embedding results of \textit{eurodist} dataset .

\subsection*{C. Several Useful Lemmas}

%\begin{Lemm}
%\label{Lemm:trace}
%Let $U, V \in \mathbb{R}^{p\times r}$, $H = U-V$. Suppose that $\mathrm{tr}(U^TV)\geq 0$, then
%\[
%\mathrm{tr}((U^TV +V^TU)H^TH) \geq 0.
%\]
%\end{Lemm}

\begin{Lemm}[\cite{Bhatia-1987}]
\label{Lemm:trace}
Let $A$ and $B$ be two positive semi-definite matrices with the size $p\times p$. Assume that $A$ is full rank, then
\[
\mathrm{tr}(AB) \geq \sigma_{\min}(A)\mathrm{tr}(B).
\]
\end{Lemm}

\begin{Lemm}
\label{Lemm:dist-twomat}
For any $U \in \mathbb{R}^{p\times r},$ let $X = UU^T$, $X_r^* = U_r^*{U_r^*}^T$, the following hold:
\begin{enumerate}
\item[(a)] \textbf{(Upper bound)} $\|X-X_r^*\|_F^2 \leq 2(\|X\|_F+\|X_r^*\|_F) \cdot \|U-U_r^*\|_F^2$, and

\item[(b)] \textbf{(Lower bound)} if $\|U-U_r^*\|_F \leq \gamma \sigma_{r}(U_r^*)$ for some $0<\gamma <1$, then
\[\|X-X_r^*\|_F^2 \geq 2(\sqrt{2}-1)\sigma_r^2(U_r^*) \|U-U_r^*\|_F^2.\]
\end{enumerate}
\end{Lemm}
%Lemma \ref{Lemm:dist-twomat}(b) is an amendment of \cite[Lemma 5.4]{Recht-2016} from the metric \footnote{where ${\cal O}$ is a set of orthogonal matrices of size $r\times r$} $\mathrm{dist}(U,V):= \min_{R\in {\cal O}} \|U-V\cdot R\|_F$ to $\|U-V\|_F$. Thus, their proofs are very similar.

\textbf{Proof.}
%\begin{proof}
\textbf{(a)} Note that
\begin{align*}
\|UU^T - U_r^*{U_r^*}^T\|_F
&= \|U(U-U_r^*)^T + (U-U_r^*){U_r^*}^T\|_F \\
&\leq (\|U\|_F + \|U_r^*\|_F)\|U-U_r^*\|_F.
\end{align*}
Thus,
$\|X-X_r^*\|_F^2
\leq (\|U\|_F + \|U_r^*\|_F)^2\|U-U_r^*\|_F^2 \leq 2(\|U\|_F^2 + \|U_r^*\|_F^2)\|U-U_r^*\|_F^2.
$

\textbf{(b)}
%By the assumption of $\|U-U_r^*\|_F \leq \gamma \sigma_{r}(U_r^*)$ with $0<\gamma <1$, we have
For any $x\in \mathbb{R}^r$, note that
\begin{align}
\label{Eq:tr-UV}
2x^T U^TU_r^* x
&= \|Ux\|_2^2 + \|U_r^* x\|_2^2 - \|(U-U_r^*)x\|_2^2 \nonumber\\
& \geq \|U_r^* x\|_2^2 - \|U-U_r^*\|_2^2 \cdot \|x\|_2^2 \nonumber\\
& \geq (1- \gamma^2) \sigma_r{(X_r^*)} \|x\|_2^2 \nonumber\\
& \geq 0 \quad (\because 0<\gamma<1),
\end{align}
where the first inequality is due to $\|Ux\|_2^2 \geq 0$ and $\|(U-U_r^*)x\|_2 \leq \|U-U_r^*\|_2 \cdot \|x\|_2$,
and the second inequality holds for $\|U_r^* x\|_2^2 \geq \sigma_r(X_r^*)\|x\|_2^2$ and $\|U-U_r^*\|_2^2 \leq \|U-U_r^*\|_F^2 \leq \gamma^2 \sigma_{r}(X_r^*)$ by the assumption of this lemma.
Thus, \eqref{Eq:tr-UV} implies
\begin{align}
\label{Eq:pd-u'ur*}
U^T U_r^* \succ 0,
\end{align}
and $U^T U_r^*$ is full rank.
%\begin{align}
%\label{Eq:tr-UV}
%\mathrm{tr}(U^TU_r^*) = \langle U, U_r^* \rangle = \frac{1}{2} (\|U\|_F^2 + \|U_r^*\|_F^2 - \|U-U_r^*\|_F^2) \geq \frac{1}{2} (\|U_r^*\|_F^2 - \gamma^2 \sigma_r^2(U_r^*)) \geq 0.
%\end{align}
Based on \eqref{Eq:pd-u'ur*}, we prove part (b).
Let $H = U-U_r^*$. Then
\begin{align*}
\|X-X_r^*\|_F^2
&=  \mathrm{tr}\left((H^TH)^2 + 4H^THH^TU_r^*\right)\\
&+  \mathrm{tr}\left(2(H^T U_r^*)^2 + 2{U_r^*}^T {U_r^*} H^TH\right).
\end{align*}
Thus, establishing Lemma \ref{Lemm:dist-twomat}(b) is equivalent to show that
\begin{align*}
&\mathrm{tr}\left((H^TH)^2 + 4H^THH^T {U_r^*} + 2(H^T {U_r^*})^2\right) \\
&+ \mathrm{tr}\left(2{U_r^*}^T {U_r^*} H^TH - c H^TH\right) \\
&\geq 0,
\end{align*}
where $c := 2(\sqrt{2}-1)\sigma_r^2({U_r^*})$. By some simple derivations, we can observe that
\begin{align*}
&\mathrm{tr}\left((H^TH)^2 + 4H^THH^T {U_r^*} + 2(H^T {U_r^*})^2\right) \\
&+ \mathrm{tr}\left(2{U_r^*}^T {U_r^*} H^TH - c H^TH\right) \\
&= \mathrm{tr}\left((H^TH+\sqrt{2}H^T {U_r^*})^2 + (4-2\sqrt{2})H^THH^T {U_r^*}\right) \\
&+ \mathrm{tr}\left(2{U_r^*}^T{U_r^*}H^TH - c H^TH\right)\\
&\geq \mathrm{tr}\left((4-2\sqrt{2})H^THH^T{U_r^*} + 2{U_r^*}^T{U_r^*}H^TH - c H^TH\right)\\
& = \mathrm{tr}\left(\left((4-2\sqrt{2})H^T{U_r^*} + 2{U_r^*}^T{U_r^*} - c {\bf I}\right)\cdot H^TH\right).
\end{align*}
Recalling $H^T {U_r^*} = U^T {U_r^*}- {U_r^*}^T{U_r^*}$, we have
\begin{align*}
&\mathrm{tr}\left(\left((4-2\sqrt{2})H^T{U_r^*} + 2{U_r^*}^T{U_r^*} - c {\bf I}\right)\cdot H^TH\right)\\
& = \mathrm{tr}\left((4-2\sqrt{2})U^T{U_r^*} \cdot H^TH\right)\\
& + \mathrm{tr}\left(\left(2(\sqrt{2}-1){U_r^*}^T{U_r^*} - c {\bf I}\right)\cdot H^TH\right)\\
& = \mathrm{tr}\left((2-\sqrt{2})(U^T{U_r^*}+{U_r^*}^TU)\cdot H^TH\right)\\
& + \mathrm{tr}\left(\left(2(\sqrt{2}-1){U_r^*}^T{U_r^*} - c {\bf I} \right)\cdot H^TH\right)\\
& \geq \mathrm{tr}\left(\left(2(\sqrt{2}-1){U_r^*}^T{U_r^*} - c {\bf I} \right)\cdot H^TH\right) \\
& \geq 0, \quad (\because c :=2(\sqrt{2}-1)\sigma_r^2({U_r^*}), \ \text{Lemma} \ \ref{Lemm:trace})
\end{align*}
where the second equality is due to $\mathrm{tr}(U^T{U_r^*} H^T H) = \mathrm{tr}(H^TH{U_r^*}^TU) = \mathrm{tr}({U_r^*}^TU H^TH)$, and the first inequality holds for \eqref{Eq:pd-u'ur*} and Lemma \ref{Lemm:trace}.
Therefore, the above inequality implies
\begin{align*}
&\mathrm{tr}\left((H^TH)^2 + 4H^THH^T{U_r^*} + 2(H^T{U_r^*})^2 \right)\\
& + \mathrm{tr}\left( 2{U_r^*}^T{U_r^*}H^TH - c H^TH\right) \geq 0,
\end{align*}
which concludes part (b) of this lemma.
%\end{proof}
$\Box$

The following lemma is similar to  \cite[Lemma 18]{Sanghavi-FGD2016}.
\begin{Lemm}
\label{Lemm:X-Xr*}
Let $X = UU^T$ and $X_r^* = U_r^*{U_r^*}^T$ be two $p\times p$ rank-$r$ positive semidefinite matrices. Let $\|U-U_r^*\|_F \leq \gamma \sigma_r(U_r^*)$ for some constant $0<\gamma <1$. Then
\[
\|X-X_r^*\|_F \leq (2\gamma + \gamma^2)\cdot \tau(U_r^*)\cdot \sigma_r(X_r^*),
\]
where $\tau(U_r^*):=\frac{\sigma_1(U_r^*)}{\sigma_r(U_r^*)}$.
\end{Lemm}

%\begin{proof}
\textbf{Proof.}
Note that
\begin{align*}
\|X-X_r^*\|_F
&= \|U(U - {U_r^*})^T + (U - U_r^*) {U_r^*}^T\|_F \\
&\leq \|U-U_r^*\|_F (\|U\|_2+\|U_r^*\|_2) \\
&\leq (2\|U_r^*\|_2+\gamma \sigma_r(U_r^*))\|U-U_r^*\|_F \\
&\leq (2+\gamma)\gamma \cdot \|U_r^*\|_2 \cdot \sigma_r(U_r^*),
\end{align*}
where the first inequality holds for the triangle inequality, Cauchy-Schwartz inequality and the fact that the spectral norm is invariant with respect to the orthogonal transformation, the second inequality is due to the following sequence of inequalities, based on the hypothesis of the lemma:
\[
\|U\|_2 - \|U_r^*\|_2 \leq \|U-U_r^*\|_2 \leq \|U-U_r^*\|_F \leq \gamma \sigma_r(U_r^*),
\]
and the last inequality holds for the fact $\sigma_r(U_r^*) \leq \|U_r^*\|_2$ and the assumption of this lemma.
The above inequality directly implies the claim of this lemma by the definition of $\tau(U_r^*)$.
%\end{proof}
$\Box$

Moreover, we need modify \cite[Lemma 19]{Sanghavi-FGD2016} as follows.
\begin{Lemm}
\label{Lemm:sigma_min_X}
Let $X = UU^T$ and $X_r^* = U_r^*{U_r^*}^T$ be two $p\times p$ rank-$r$ positive semidefinite matrices. Let $\|U-U_r^*\|_F \leq \gamma \sigma_r(U_r^*)$ for some constant $0<\gamma <1$. Then
\[
\sigma_r(U) \geq (1-\gamma) \sigma_r(U_r^*).
\]
\end{Lemm}

\textbf{Proof.}
%\begin{proof}
Using the norm ordering $\|\cdot\|_2 \leq \|\cdot\|_F$ and the Weyl's inequality for perturbation of singular values (see, \cite[Theorem 3.3.16]{Horn-matrixbook1991}), we get
\[
|\sigma_i(U)-\sigma_i(U_r^*)| \leq \gamma \sigma_r(U_r^*), \ 1\leq i \leq r,
\]
which implies that
\[
\sigma_r(U) \geq (1-\gamma)\sigma_r(U_r^*).
\]
%Thus, $\sigma_r(X) \geq (1-\sqrt{\gamma})^2 \sigma_r(X_r^*).$
%\end{proof}
$\Box$

\begin{Lemm}
\label{Lemm:feasibility}
Let Assumption \ref{Assump:objfun} hold. Let $X = UU^T$ and $X_r^* = U_r^*{U_r^*}^T$ be two $p\times p$ rank-$r$ ($r<p$) positive semidefinite matrices. Suppse that $\|U-U_r^*\|^2_F < \gamma_0 \sigma_r(X_r^*)$, where $\gamma_0$ is specified in \eqref{Eq:kappa}. Then the following hold:
\begin{enumerate}
\item[(a)] (\textbf{Bounded gradient})
\begin{align*}
&\|\nabla f(X)\|_F \\
&\leq  \|\nabla f(X_r^*)\|_F + (2\sqrt{\gamma_0}+\gamma_0)L\tau(U_r^*)\sigma_r(X_r^*),
\end{align*}

\item[(b)] (\textbf{Feasibility of} $\bar{X}$) Let $\bar{X}:=X - \frac{\bar{\eta}}{L} {\cal P}_U \nabla f(X) {\cal P}_U $, where $\bar{\eta}$ is specified in \eqref{bar eta}, then $\bar{X}$ is a feasible point, i.e., $\bar{X}$ is symmetric and positive semidefinite.

\item[(c)] $ ({\bf I}-{\cal P}_U)X_r^* =0.$
\end{enumerate}
\end{Lemm}
%\begin{proof}
\textbf{Proof.}
\textbf{(a)} Note that
\begin{align*}
\|\nabla f(X)\|_F
&\leq \|\nabla f(X_r^*)\|_F + L\|X-X_r^*\|_F\\
&\leq \|\nabla f(X_r^*)\|_F + (2\sqrt{\gamma_0}+\gamma_0)L\tau(U_r^*)\sigma_r(X_r^*),
\end{align*}
where the first inequality holds for the $L$-Lipschitz differentiability of $f$, and the second inequality holds for Lemma \ref{Lemm:X-Xr*}.

\textbf{(b)} Since ${\cal P}_U X {\cal P}_U = X$ and $X$ is rank-$r$, then
\[X - \frac{\bar{\eta}}{L} \cdot {\cal P}_U \nabla f(X){\cal P}_U  = {\cal P}_U(X - \frac{\bar{\eta}}{L} \cdot \nabla f(X)){\cal P}_U,\]
which implies that $\bar{X}$ is symmetric and that the last $p-r$ eigenvalues of  the matrix $X - \frac{\bar{\eta}}{L} \cdot {\cal P}_U \nabla f(X){\cal P}_U$ are zero, that is, $\lambda_i(X - \frac{\bar{\eta}}{L} \cdot {\cal P}_U \nabla f(X){\cal P}_U ) =0$ for $i=r+1,\ldots,p$. While for any $i=1,\ldots, r$,
\begin{align*}
&\lambda_i\left(X - \frac{\bar{\eta}}{L} \cdot {\cal P}_U \nabla f(X){\cal P}_U \right)\\
& \geq \lambda_i(X) - \frac{\bar{\eta}}{L} \cdot \lambda_{\max}({\cal P}_U \nabla f(X){\cal P}_U )\\
& \geq \sigma_r(X) - \frac{\bar{\eta}}{L} \cdot \sigma_{\max}(\nabla f(X))\\
& \geq (1-\sqrt{\gamma_0})^2 \sigma_r(X_r^*) \\
& - \frac{\bar{\eta}}{L} \cdot \left[\|\nabla f(X_r^*)\|_F + (2\sqrt{\gamma_0}+\gamma_0)L\tau(U_r^*)\sigma_r(X_r^*)\right] \\
&\geq 0,
\end{align*}
where the third inequality holds for Lemma \ref{Lemm:sigma_min_X} and \textbf{(a)} of this lemma, and the final inequality holds for the definition of $\bar{\eta}$. Therefore, $\bar{X}$ is positive semidefinite.

\textbf{(c)} By $\|U-U_r^*\|_F < \sqrt{\gamma_0}\sigma_r(U_r^*)$ and $0<\sqrt{\gamma_0}<1$, we have
\[
\sigma_i(U) \cdot \sigma_{i}(U_r^*) > 0, \quad i\in \{1,\ldots,r\},
\]
and
\[
\sigma_{i}(U_r^*)=0, \ \sigma_{i}(U)=0, \quad i \in \{r+1,\ldots,p\},
\]
which implies that $U_r^*$ lies in the subspace spanned by $U$. In other words, $U_r^*$ does not lie in the orthogonal subspace of the subspace spanned by $U$, that is, the following holds
\[
({\bf I}-{\cal P}_U)U_r^* =0.
\]
Thus, $({\bf I}-{\cal P}_U)X_r^* =0$.
%\end{proof}
$\Box$

\end{document}